\documentclass[12pt]{amsart}
\pdfoutput=1
\usepackage[margin=1in]{geometry}
\usepackage{graphicx} % Enhanced support for graphics
\usepackage[usenames]{color} % Colour control for LATEX documents
\usepackage{booktabs} % Publication quality tables
\usepackage{dcolumn} % Align on the decimal point of numbers in tabular columns
\usepackage{enumitem} % Control layout of itemize, enumerate, description
\usepackage[colorlinks=true, linkcolor=blue, citecolor=magenta]{hyperref} % Extensive support for hypertext
\usepackage{cite} % Improved citation handling in LATEX

\usepackage{mathtools} % An extension package to amsmath (also correct various bugs/defeciencies)
\usepackage{amssymb} % Also loads amsfonts
\usepackage{mathrsfs} % Script math fonts
\usepackage{bbm} % Capital letters with double-striked left side and normal right side
\usepackage{amsthm} %
\usepackage{caption}
\usepackage{caption}

\theoremstyle{definition}

\theoremstyle{remark}

\numberwithin{equation}{section}
\usepackage[T1]{fontenc}
\usepackage{microtype} % Subliminal refinements towards typographical perfection

\newcommand{\R}{{\mathbb R}}

\newcommand{\be}{\begin{eqnarray}}
\newcommand{\ben}{\begin{eqnarray*}}
\newcommand{\en}{\end{eqnarray}}
\newcommand{\enn}{\end{eqnarray*}}

\newcommand{\pa}{\partial}

\newcommand{\hth}{\hat{\theta}}
\newcommand{\hx}{\hat{x}}

\usepackage{subfigure}
\title[A machine learning based inverse scattering scheme]{Machine learning based data retrieval for inverse scattering problems with incomplete data}

\author{Yu Gao}
\address{School of Mathematics, Jilin University, Changchun, China}
\email{gaoyu19@mails.jlu.edu.cn}

\author{Kai Zhang$^*$}
\address{School of Mathematics, Jilin University, Changchun, China}
\email{zhangkaimath@jlu.edu.cn}

\begin{document}

\maketitle

\begin{abstract}

We are concerned with the inverse scattering problems associated
with incomplete measurement data. It is a challenging topic of
increasing importance that arise in many practical applications. Based on a
prototypical working model, we propose a machine learning based
inverse scattering scheme, which integrates a CNN (convolution
neural network) for the data retrieval. The proposed method can
effectively cope with the reconstruction under limited-aperture
and/or phaseless far-field data. Numerical experiments verify the
promising features of our new scheme.

\medskip



\end{abstract}

%%%%%%%%%%%%%%%%%%%%%%%%%%%%%%%%%%%%%%%%%%%%%%%%%%%%%%%%%%%%%%%%%%

\section{Introduction}

%%%%%%%%%%%%%%%%%%%%%%%%%%%%%%%%%%%%%%%%%%%%%%%%%%%%%%%%%%%%%%%%%%

\subsection{Background and motivation}

Inverse scattering problems are concerned with the recovery of
unknown/inaccessible objects from the corresponding wave probing
data. They lie at the heart of many scientific and industrial
developments including radar and sonar, medical imaging, geophysical
exploration and non-destructive testing
\cite{AKn,BK,CCHn,ColtonCoyleMonk,CK,CKn}.

One class of inverse
scattering problems that has received significant attentions
recently in the literature is the inverse problems with incomplete
data. Those problems could arise in various applications of great
importance. The incomplete measurement data may be of limited
aperture or/and without phase information. Various mathematical
strategies have been proposed for retrieving the full wave data from
the measured partial/incomplete data, and then they are integrated
into the reconstruction process
\cite{AmmariGarnierSolna2013,BK,DZGn,JLZn,KKNn,LLn,LLWn,SKNn,ZG,ZGLLn}.
Generally speaking, if sufficient a-priori information is
available on the underlying target scatterer, one may still be
possible to achieve the recovery to a certain content level. Otherwise, the reconstruction shall suffer from great deterioration
due to the intrinsic lack of information. However, some newly
emerging applications may have more restrictive and higher quality
requirements on the data retrieval and ultimately the recovery
effects. In fact, it is noted that there is an increasing trend on
designing portable and handheld medical devices that are fit for
family use, say e.g. ultrasound or MRI scanners that can produce
medical imaging on one's smart phone. Clearly, in addition to the
hardware requirements, the core to the success of such conceptual
applications would be the development of novel inverse scattering
schemes that can work effectively and efficiently with highly
incomplete measurement data. This motivates the current study by
integrating the machine learning techniques into the data retrieval.

Machine learning methods, in particular the convolutional neural
networks (CNNs), have been regarded as a revolutionary idea for many
applied sciences. In fact, they have been successfully developed for
dealing with various inverse problems with an increasing amount of
literature. In this article, we propose a CNN-based data retrieval
approach, which in combination with a sampling-type imaging scheme
can produce stunning reconstructions for inverse scattering problems
with limited-aperture or/and phaseless measurement data. Since it is
rather impractical to consider too many different inverse scattering
problems, we develop our method based on a prototypical model, the
so-called inverse acoustic obstacle problem. The rest of this
section is devoted to the introduction of this model inverse
scattering problem as well as some relevant discussions on the
reconstruction methods.

%%%%%%%%%%%%%%%%%%%%%%%%%%%%%%%%%%%%%%%%%%%%%%%%%%%%%%%%%%%%%%%%%%

\subsection{Mathematical setup and relevant discussions}

Let $k=\omega/c\in\mathbb{R}_+$ be the wave number of a time
harmonic acoustic wave with $\omega\in\mathbb{R}_+$ and
$c\in\mathbb{R}_+$ denoting, respectively, the angular frequency and
the wave speed in a homogeneous background space. Let
$\Omega\subset\mathbb{R}^n (n=2,\, 3)$ be a bounded domain with a
Lipschitz-boundary $\partial\Omega$ such that its complement
$\mathbb{R}^n\setminus\overline{\Omega}$ is connected. Let the
incident field $u^i$ be a plane wave of the form
\begin{equation*}
u^i: =\ u^i(x,\hat{\theta},k) = e^{ikx\cdot \hat{\theta}},\quad x\in\mathbb{R}^n\,,
\end{equation*}
where $\hat{\theta}\in \mathbb{S}^{n-1}$ denotes the impinging
direction of the incident wave and
$\mathbb{S}^{n-1}:=\{x\in\mathbb{R}^n:|x|=1\}$ is the unit sphere in
$\mathbb{R}^n$. In acoustic probing, one sends an incident wave
$u^i$ to detect the unknown/inaccessible scatterer $\Omega$. The
scatterer interrupts the propagation of the plane
wave and generates the so-called scattered wave $u^s$. Set
$u=u^i+u^s$ to denote the total wave field. The acoustic obstacle
scattering is described by the following boundary value problem,
\begin{equation}
\label{HemEquobstacle}
\begin{cases}
\ \ \Delta u + k^2 u = 0 & \hspace*{-2cm}\mbox{in}\ \ \mathbb{R}^n\setminus\overline{\Omega},\medskip\\
\ \  u = 0 & \hspace*{-2cm}\mbox{on }\ \ \partial\Omega,\medskip\\
\displaystyle{\lim_{r\rightarrow\infty}r^{\frac{n-1}{2}}\left(\frac{\partial u^{s}}{\partial r}-iku^{s}\right) =\,0,\ \ r=|x|\ \ \mbox{for}\ x\in\mathbb{R}^n.}
\end{cases}
\end{equation}
In \eqref{HemEquobstacle}, the homogeneous Dirichlet boundary
condition on $\partial D$ signifies that $D$ is a sound-soft
obstacle. The third equation is known as the Sommerfeld radiation
condition, which characterizes the outgoing nature of the scattered
wave field. The forward scattering problem is to compute $u^s$ (or
$u$) for a given incident field $u^i$ and domain $D$. We refer to
\cite{CK, Mclean} for the well-posedness of the forward scattering
problem \eqref{HemEquobstacle} in $H^1_{\rm
loc}(\mathbb{R}^n\setminus\overline{D})$.

In particular, one has the following asymptotic expansion of the scattered field,
\begin{equation}
\label{0asyrep}
u^s(\Omega; x, \theta)
=\frac{e^{i\frac{\pi}{4}}}{\sqrt{8k\pi}}\left(e^{-i\frac{\pi}{4}}\sqrt{\frac{k}{2\pi}}\right)^{n-2}\frac{e^{ikr}}{r^{\frac{n-1}{2}}}
   \left\{u^{\infty}(\Omega; \hat{x}, \hat{\theta})+\mathcal{O}\left(\frac{1}{r}\right)\right\}\quad\mbox{as }\,r\rightarrow\infty,
\end{equation}
which holds uniformly with respect to all observation directions
$\hat{x}:=x/|x|\in \mathbb{S}^{n-1}$. The complex-valued function
$u^\infty$ in \eqref{0asyrep} defined on the unit sphere
$\mathbb{S}^{n-1}$ is known as the scattering amplitude or far-field
pattern. The inverse scattering problem we are concerned with is to
recover $\Omega$ by the knowledge of the far-field pattern
\begin{equation*}
  u^{\infty}(\Omega; \hat{x}, \hat{\theta}): (\hat{x}, \hat{\theta})\in\Gamma \times\Sigma,
\end{equation*}
where $\Gamma$ and $\Sigma$ are open subsets of $\mathbb{S}^{n-1}$
and are referred to as the observation aperture and incident aperture,
respectively. If $\Gamma=\mathbb{S}^{n-1}$ and
$\Sigma=\mathbb{S}^{n-1}$, then the corresponding inverse scattering
problem is said to have full-aperture measurement data. Otherwise,
it is said to have limited-aperture measurement data. By introducing
an abstract operator $\mathcal{F}$ which sends the underlying
obstacle to its corresponding far-field pattern and is defined by
the Helmholtz system \eqref{HemEquobstacle}, the inverse problem can
be formulated as the following operator equation,
\begin{equation}\label{eq:ip1}
\mathcal{F}(\Omega)=u^\infty(\hat x, \hat\theta), \ \ (\hat x, \hat\theta)\in\Gamma\times\Sigma.
\end{equation}
It can be verified that the inverse problem \eqref{eq:ip1} is
nonlinear. Moreover, it is severely ill-conditioned in the sense of
Hadamard; that is, a small perturbation in the far-field data of
\eqref{eq:ip1} may cause a significant change in reconstructing
$\Omega$ (cf. \cite{DR,LPRX,R1}).

It is known that the full-aperture data uniquely determines the
obstacle $D$ (cf. \cite{CK}). Since $u_\infty(\hat{x},
\hat{\theta})$ is a real-analytic function on
$\mathbb{S}^{n-1}\times\mathbb{S}^{n-1}$ (cf. \cite{CK}), the
far-field pattern on $\Gamma \times \Sigma$ can be extended to the
full-aperture data by unique continuation. Hence the
limited-aperture data also uniquely determines the obstacle.
However, it is well known that the analytic continuation is a
severely ill-conditioned process (cf. \cite{Atkinson}). Thus the
inverse scattering problem suffers from an increasing level of
ill-posedness as the size of aperture decreases. Indeed, this
phenomenon has been observed and investigated in the existing
literature, see e.g. \cite{AhnJeonMaPark, AmmariGarnierSolna2013,
BaoLiu, ColtonCoyleMonk, IkehataNiemiSiltanen, L4,
LiWangWangZhao, MagerBleistein,MagerBleistein1978,Robert1987,
Zinn1989}. In what follows we shall also present several numerical
examples to reemphasize this point.

There is another scenario of practical importance where the measurement data is given by
\begin{equation}\label{eq:phaseless1}
\mathcal{F}(\Omega)= \left|u^{\infty}(\Omega; \hat{x}, \hat{\theta})\right|,\ \ (\hat{x}, \hat{\theta})\in\Gamma \times\Sigma.
\end{equation}
For the phaseless inverse problem \eqref{eq:phaseless1}, there is a
well-known obstruction that one cannot determine the location of the
obstacle $\Omega$. That is, if one shifts the obstacle to another
location, the modulus of the corresponding far-field pattern remains
unchanged. However, it is still possible to determine only the shape
of the obstacle by knowing its position a priori (cf.
\cite{LLn,LLWn}). Many inverse scattering schemes developed for
phased measurement data in principle do not work for the phaseless
case. Hence, the phaseless inverse scattering problem constitutes
another challenging topic in the literature
\cite{DZGn,JLZn,JLZn2,LLn,LLWn,ZG,ZGLLn}. As indicated in
\cite{ZG,ZGLLn}, in the phaseless setup, it is physically
more relevant to consider the measurement of the modulus of the
total field on a sufficient large surface enclosing the obstacle
other than the modulus of the far-field pattern. We would like to
emphasize that the data retrieval method developed in what follows
can be easily extended to work for the case with the modulus of the
total wave field. However, in order to have a uniform formulation as
well as to ease our exposition, we shall confine our study in the
phaseless case with the measurement data given by
\eqref{eq:phaseless1}.

The rest of the paper is organized as follows. In Section 2, we
consider a sampling-type method for the inverse obstacle problem and
show by a few numerical examples the intrinsic ill-posedness due to
the lack of information. In Section 3, we briefly go through the
main ingredients on CNN for our subsequent use. Section 4 is devoted
to the algorithmic developments and numerical experiments.

%%%%%%%%%%%%%%%%%%%%%%%%%%%%%%%%%%%%%%%%%%%%%%%%%%%%%%%%%%%%%%%%%%

\section{A sampling-type reconstruction scheme and its localised features}

The mainstream numerical methods for solving the inverse scattering
problems \eqref{eq:ip1} or \eqref{eq:phaseless1} can be classified
as sampling-type methods and optimization based iterative methods.
Most of the sampling methods, e.g. MUSIC-type methods
\cite{AKn,AILPn,AILn}, linear sampling method \cite{CCMn} and
factorization method \cite{KirschGrinberg}, are designed to work
with the full-aperture data. Some of them can be modified to cope
with limited-aperture data, but only under special circumstances and
suffer from the intrinsic ill-posedness. The optimization methods
can be applied for any limited-aperture data as long as the
objective functional is properly defined. The limited-aperture case
has been studied by researchers using other methods, such as a regularized homotopy continuation
method \cite{BaoLiu}, a variant of the enclosure method \cite{IkehataNiemiSiltanen}, a generalization of the orthogonal
projection method \cite{Robert1987}, and the methods based on transformed field
expansion \cite{LiWangWangZhao}. The methods for reconstructing polygonal and
polyhedral obstacles also work with limited-aperture or even
phaseless data \cite{LLn,LLWn}.

Another approach for solving inverse scattering problems with
limited aperture data is to first retrieve the missing data and then
apply the existing numerical methods for the corresponding
full-aperture problems. In the literature, the missing information
retrieval is mainly based on the use of the intrinsic and generic
structures of the far-field data. The aforesaid structures include
the real-analyticity and reciprocity relations that hold generically
for any set of far-field data, see e.g. \cite{LiuSun17}.
There are also some strategies developed based on introducing an
artificial reference scatterer into the scattering system for
retrieving the phase information of the measurement data in the
phaseless case \cite{DZGn,JLZn,JLZn2,LLZn2,ZG,ZGLLn}.

\begin{figure}[!htbp]
\centering
\includegraphics[width=1.8in]{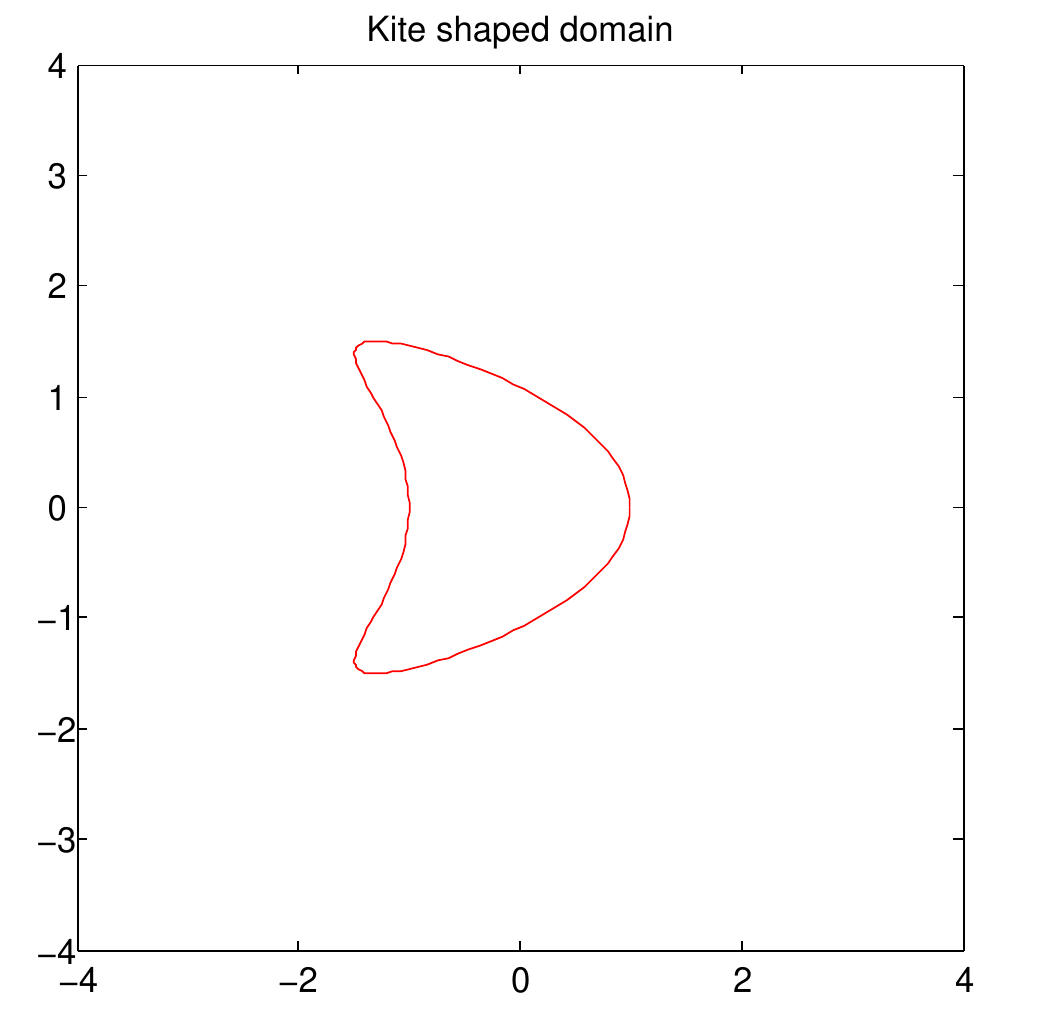}\qquad\qquad
\includegraphics[width=1.8in]{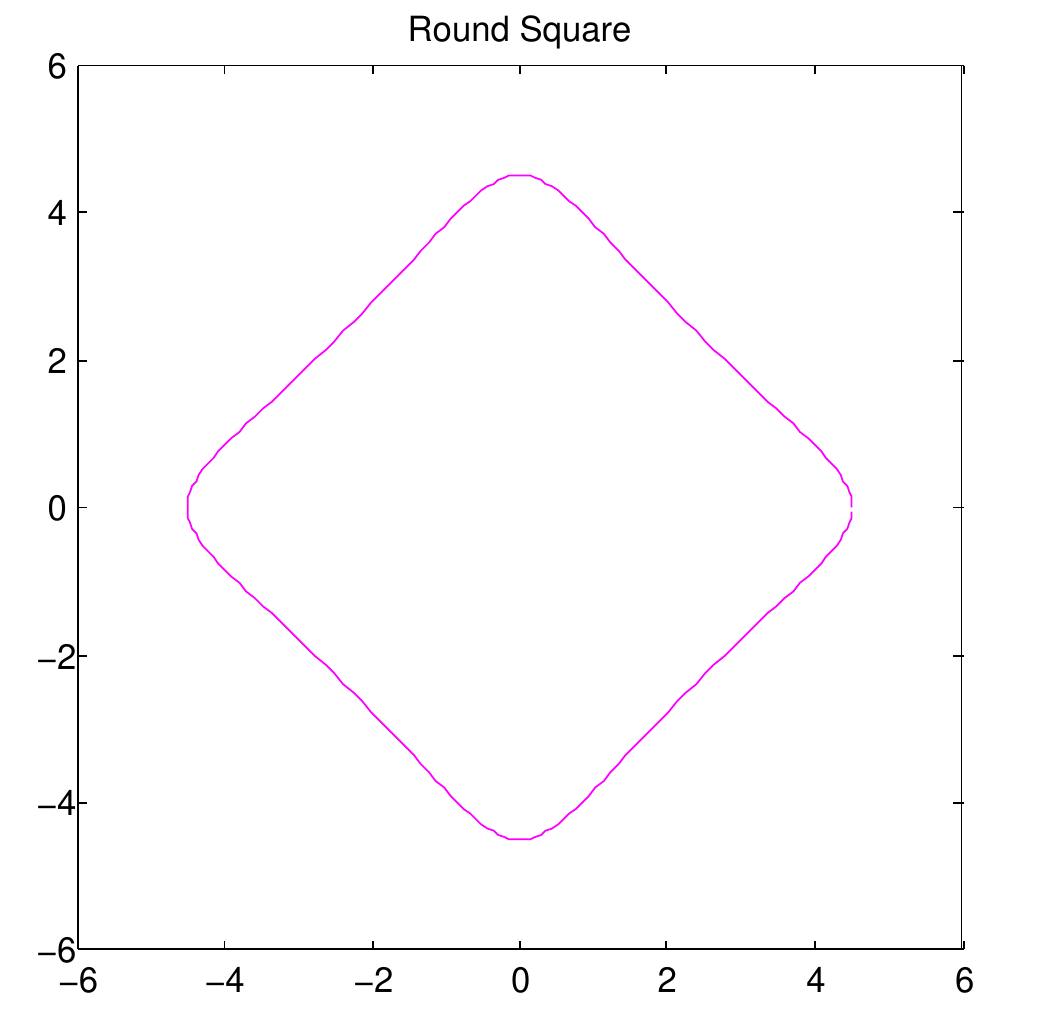}
\caption{True obstacles: kite shaped obstacle and round-square shaped obstacle.}
\label{truedomains}
\end{figure}

In what follows, we present a sampling-type method for the inverse
scattering problem \eqref{eq:ip1}, and use several examples to
discuss the intrinsic ill-posedness of \eqref{eq:ip1} in the case
with missing information. The core of the aforesaid method is the
following imaging functional

\be\label{IndicatorDSM} I(z):=\Big|\int_{\Sigma}e^{-ik\hth\cdot
z}\int_{\Gamma}u^\infty(\hx,\hth)e^{ik\hx\cdot
z}ds(\hx)ds(\hth)\Big|^{2},\quad z\in\R^n. \en The characteristic
behaviour of the imaging functional \eqref{IndicatorDSM} has been
studied in \cite{N1, LLWn1,LLZn, LiuIP17, LiuSun17}, mainly in the case with
full measurement data. Next, we present some numerical examples of
the obstacle reconstruction by using the imaging functional
\eqref{IndicatorDSM} in the two dimensional case. The boundaries of
the obstacles used in our numerical experiments are parameterized as
follows: \ben
\label{kite}&\mbox{\rm kite:}&\quad x(t)\ =(a,b)+\ (\cos t+0.65\cos 2t-0.65, 1.5\sin t),\quad 0\leq t\leq2\pi,\qquad\\
%\label{peanut}&\mbox{\rm Peanut:}&\quad x(t)\ =(a,b)+\, \sqrt{3\cos^2 t+1}(\cos t, \sin t),\quad 0\leq t\leq2\pi,\\
%\label{pear}&\mbox{\rm Pear:}&\quad x(t)\ =(a,b)+(2+0.3\cos 3t)\, (\cos t, \sin t),\quad 0\leq t\leq2\pi,\\
\label{circle}&\mbox{\rm round square:}&\quad x(t)\ =(a,b)+9/4
(\cos^{3} t+\cos t, \sin^{3} t+\sin t),\quad 0\leq t\leq2\pi, \enn
with $(a,b)$ signifying the center of the obstacle which may be
different in different examples. The plot of the regions centered at
$(0,0)$ is shown in Figure \ref{truedomains}.

\begin{figure}[!htbp]
\centering
\hspace*{-1cm}\subfigure[\textbf{$\phi\in (0,\pi/2)$}]{
    \includegraphics[width=2in]{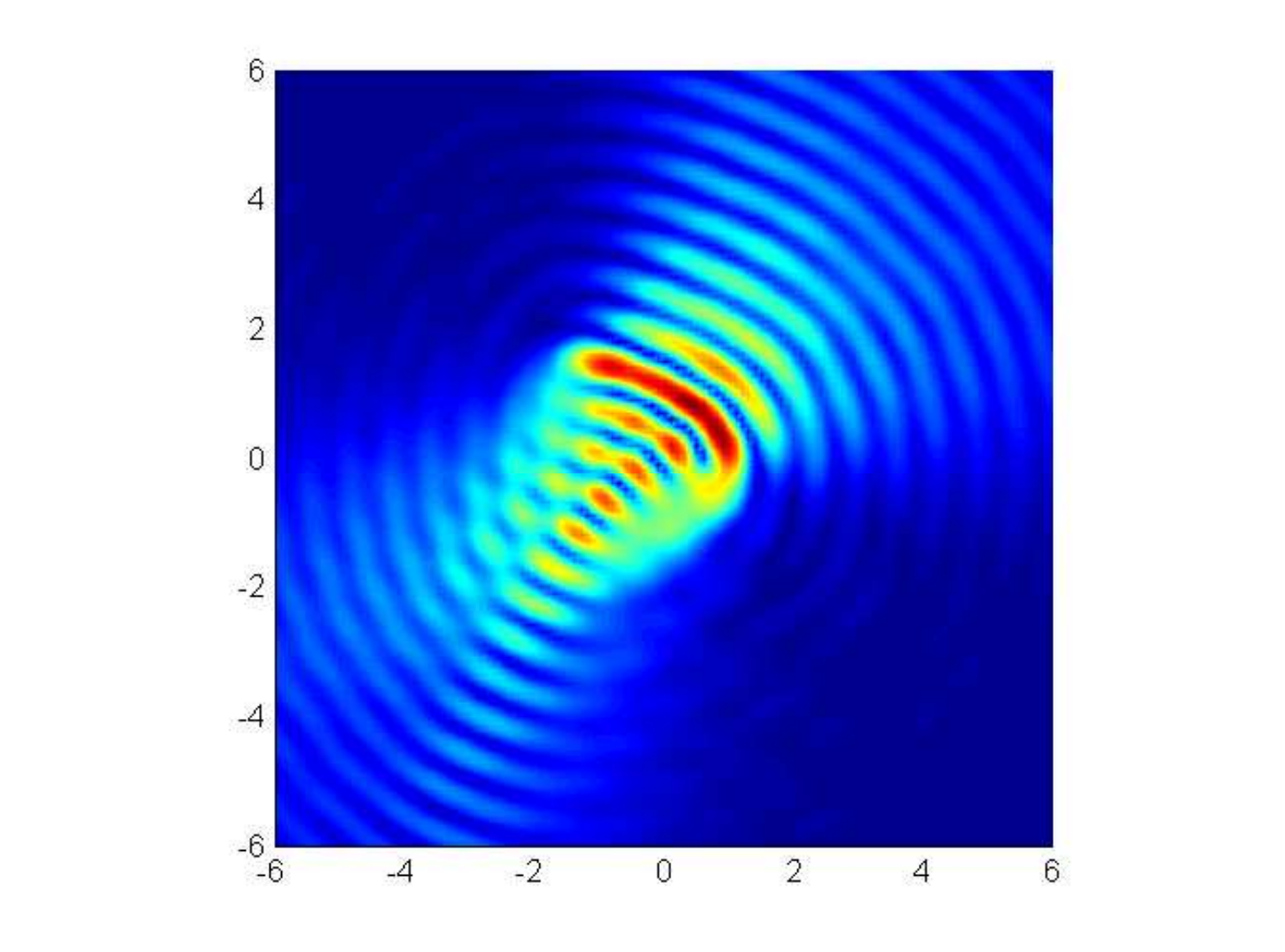}}\hspace*{-1cm}
  \subfigure[\textbf{$\phi\in (\pi/2,\pi)$}]{
    \includegraphics[width=2in]{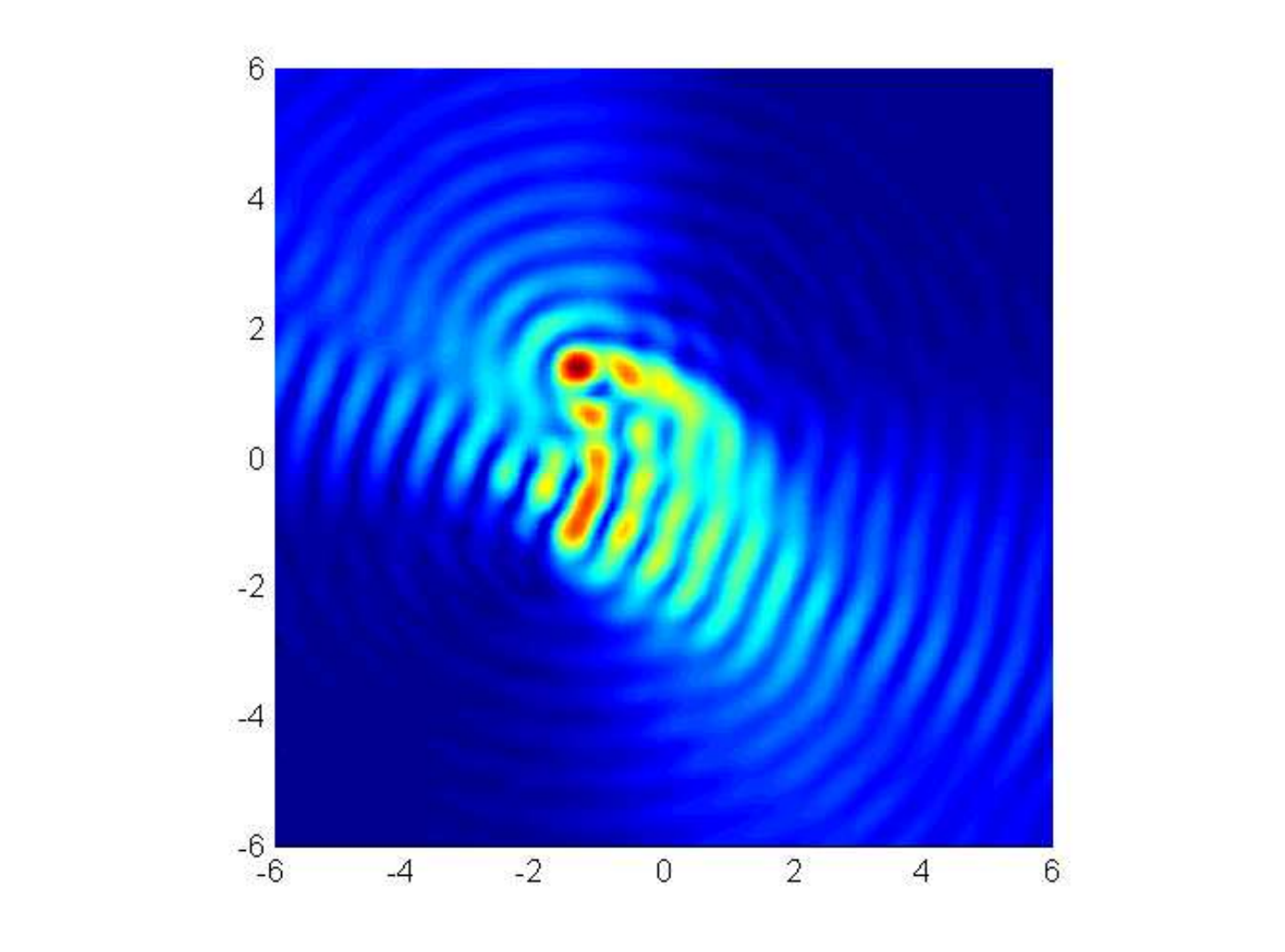}}\hspace*{-1cm}
  \subfigure[\textbf{$\phi\in (\pi,3\pi/2)$}]{
    \includegraphics[width=2in]{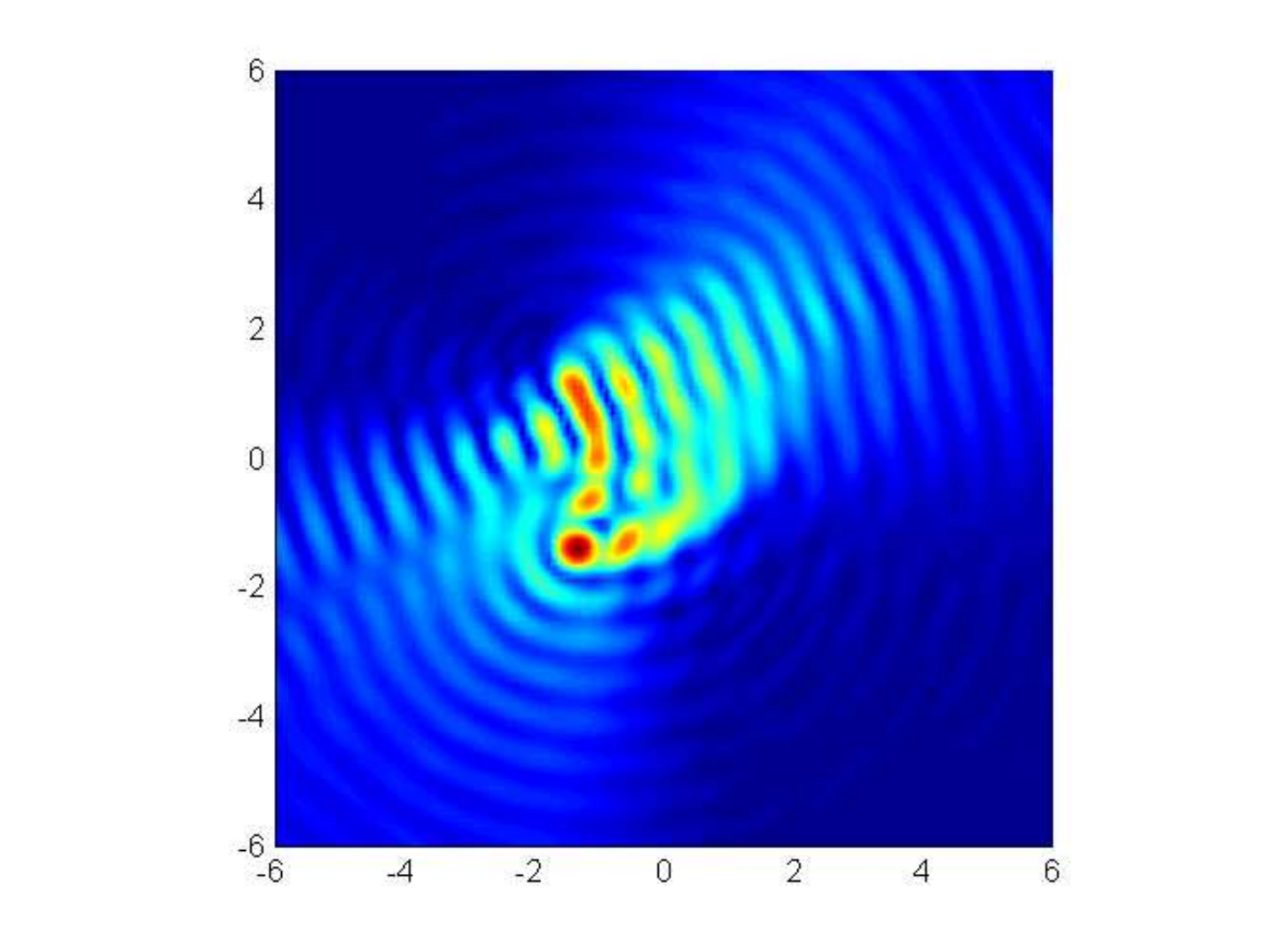}}\hspace*{-1cm}
  \subfigure[\textbf{$\phi\in (3\pi/2,2\pi)$}]{
    \includegraphics[width=2in]{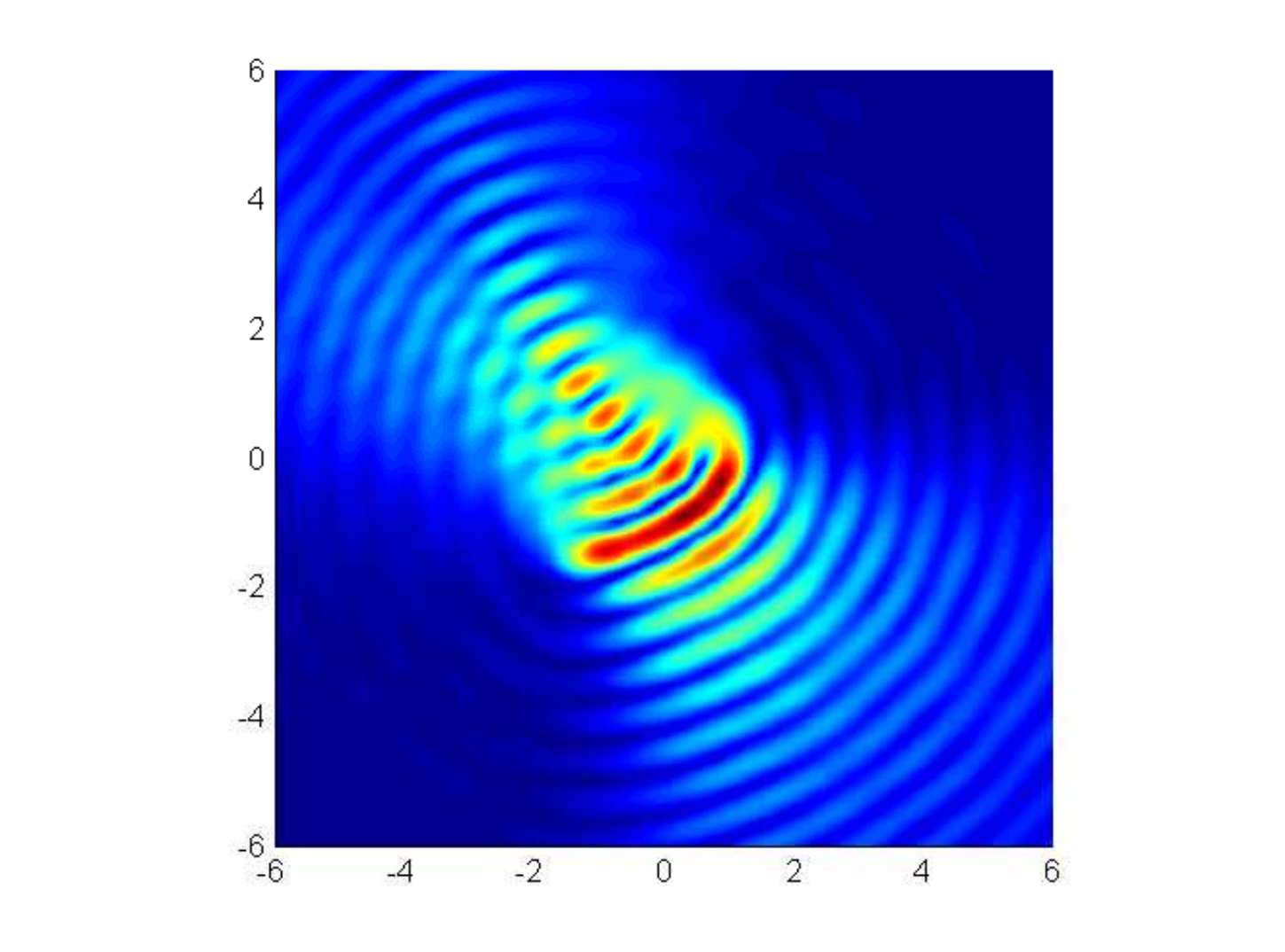}}
\caption{Reconstructions of the kite shaped sound-soft obstacle with the wave number $k=5$, $10\%$ noise and different observation apertures.}
\label{kiteapertures}
\end{figure}

\begin{figure}[!htbp]
  \centering
  \hspace*{-1cm}\subfigure[\textbf{$\phi\in (0,\pi)$}]{
    \includegraphics[width=2in]{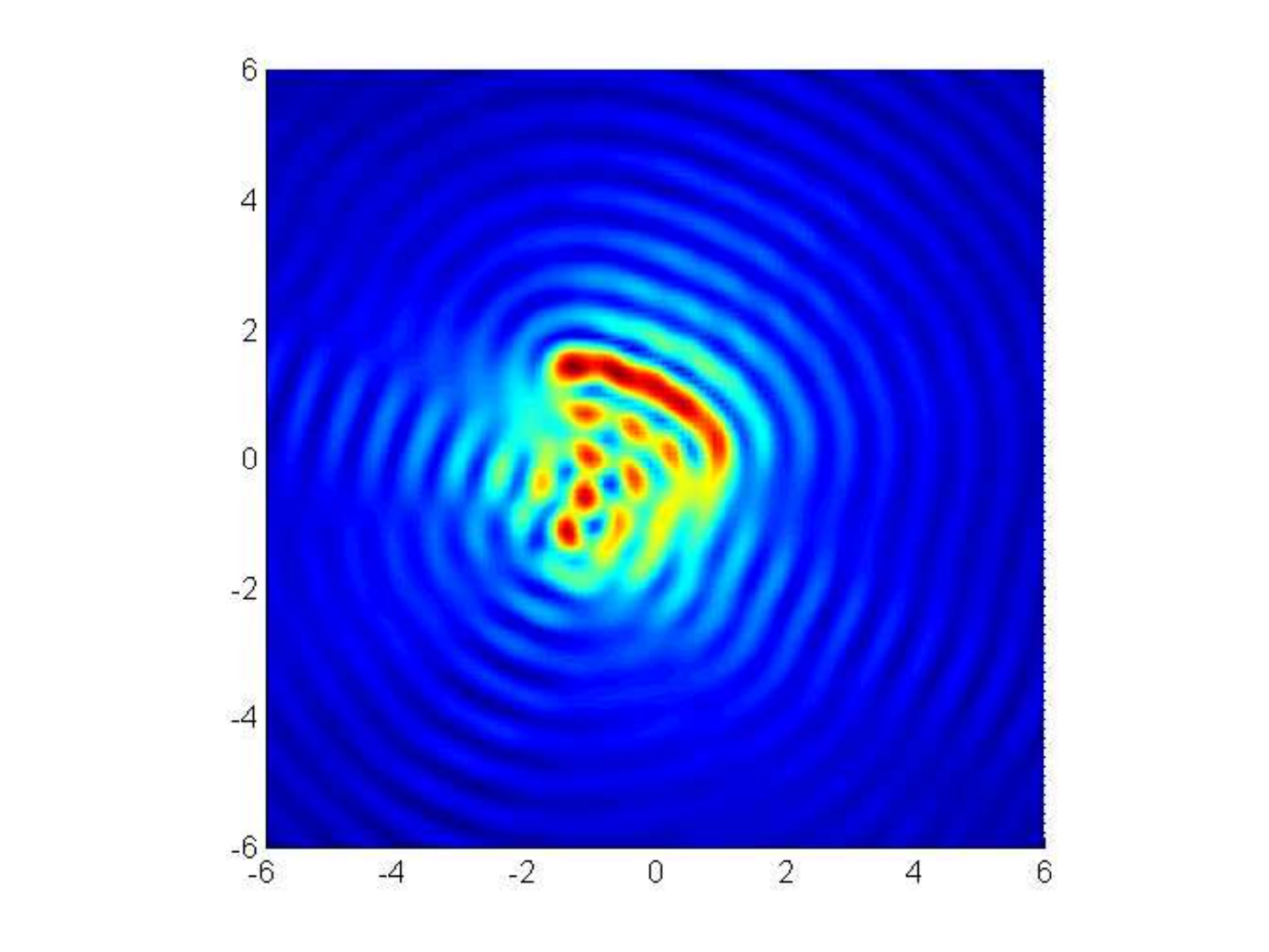}}\hspace*{-1cm}
  \subfigure[\textbf{$\phi\in (0,3\pi/4)$}]{
    \includegraphics[width=2in]{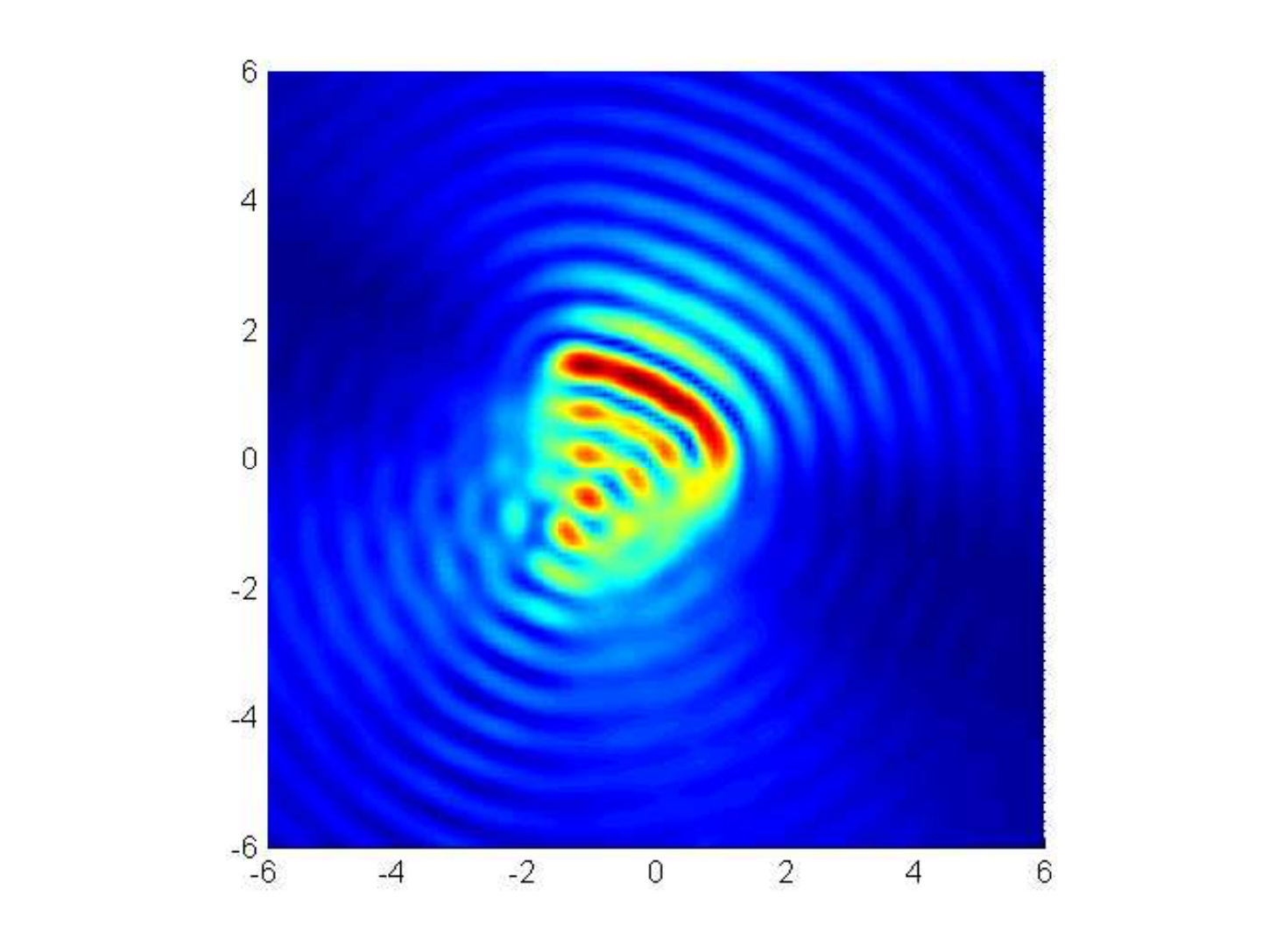}}\hspace*{-1cm}
  \subfigure[\textbf{$\phi\in (0,\pi/2)$}]{
    \includegraphics[width=2in]{pic/Kite1.pdf}}\hspace*{-1cm}
  \subfigure[\textbf{$\phi\in (0,\pi/4)$}]{
    \includegraphics[width=2in]{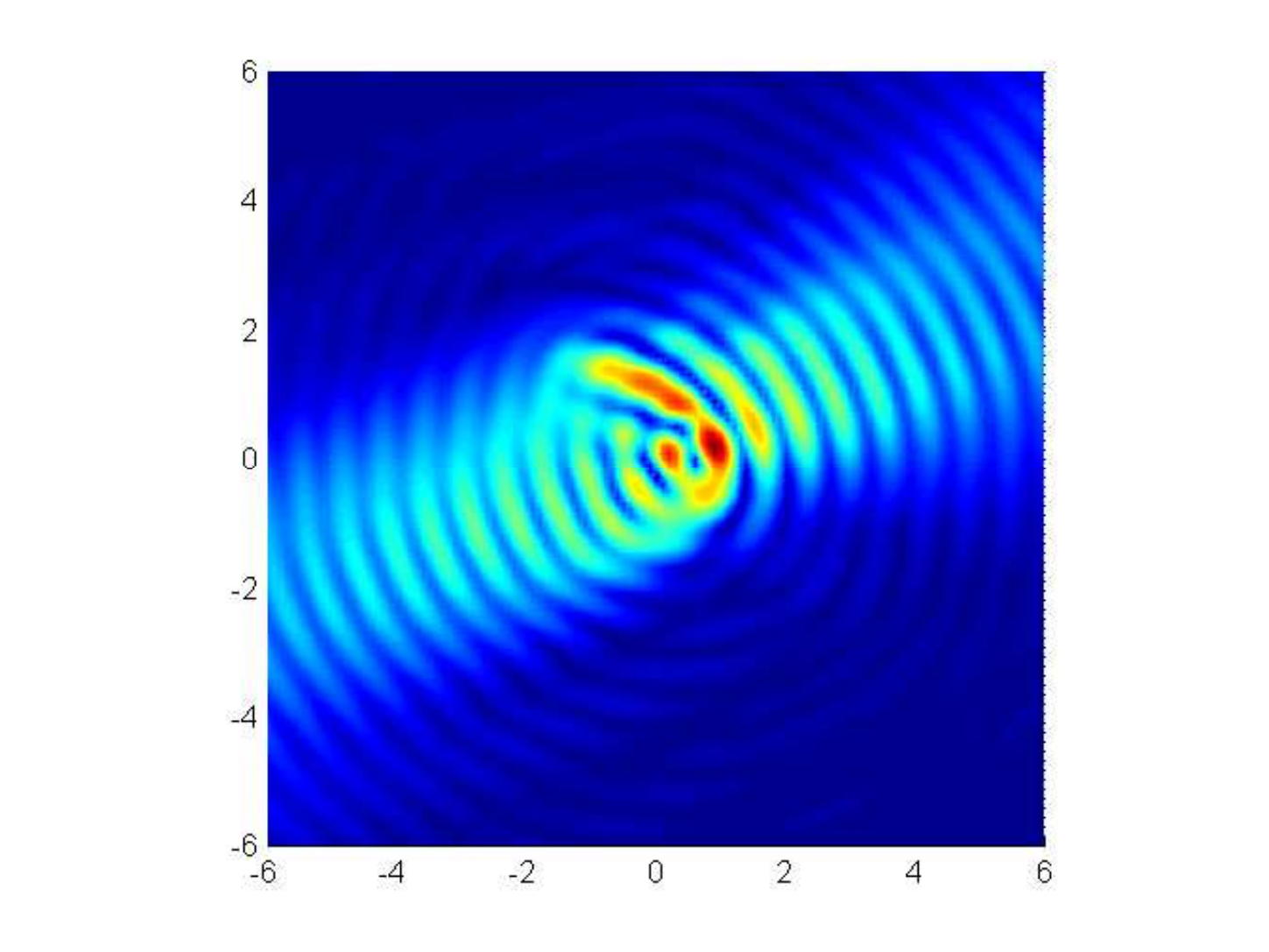}}
\caption{Reconstructions of the kite shaped obstacle with the wave number $k=5$, $10\%$ noise and different observation apertures.}
\label{kiteapertures2}
\end{figure}

\begin{figure}[!htbp]
  \centering
\hspace*{-1cm}\subfigure[\textbf{$k=3$}]{
    \includegraphics[width=2in]{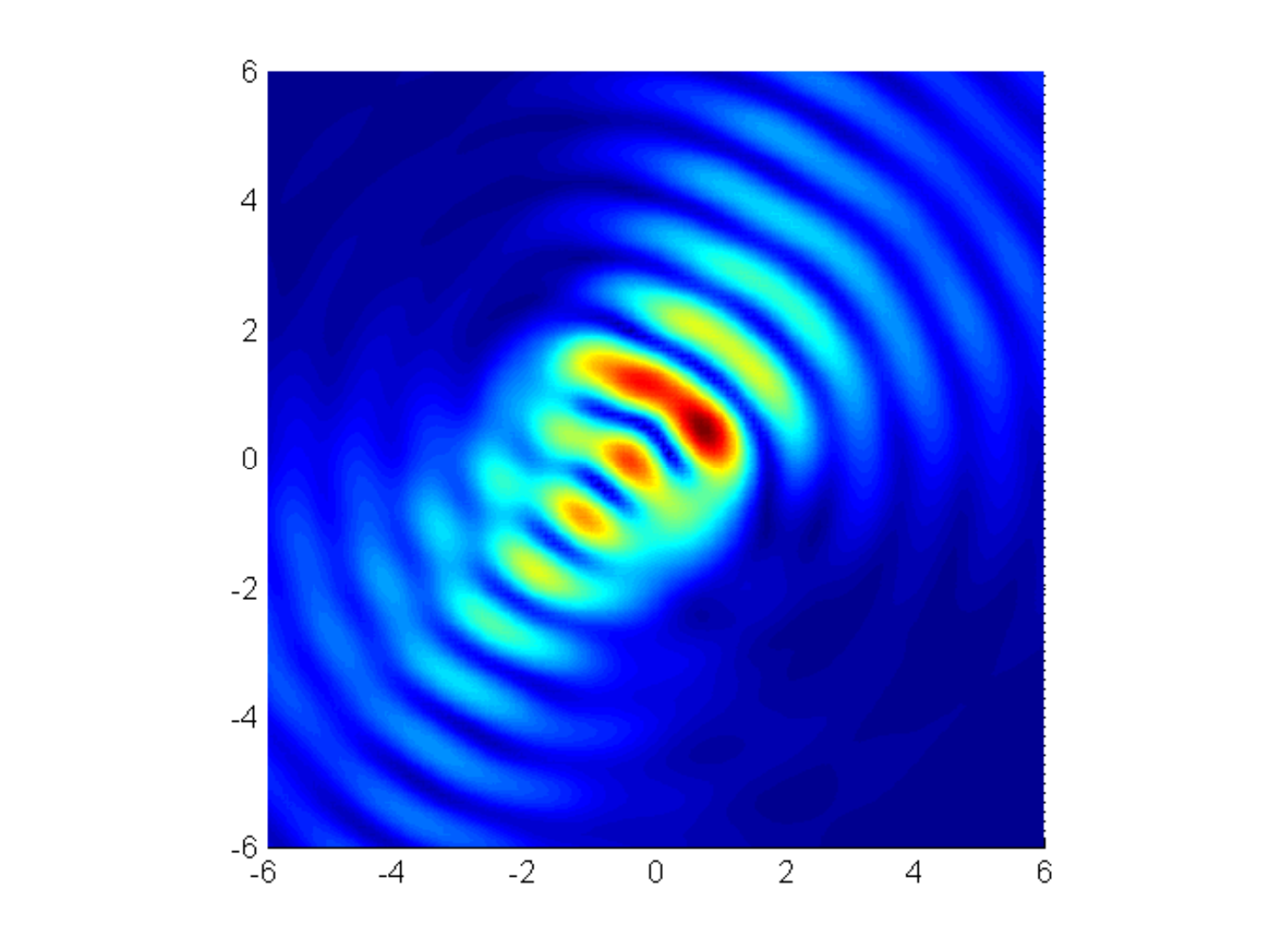}}\hspace*{-1cm}
  \subfigure[\textbf{$k=5$}]{
    \includegraphics[width=2in]{pic/Kite1.pdf}}\hspace*{-1cm}
  \subfigure[\textbf{$k=7$}]{
    \includegraphics[width=2in]{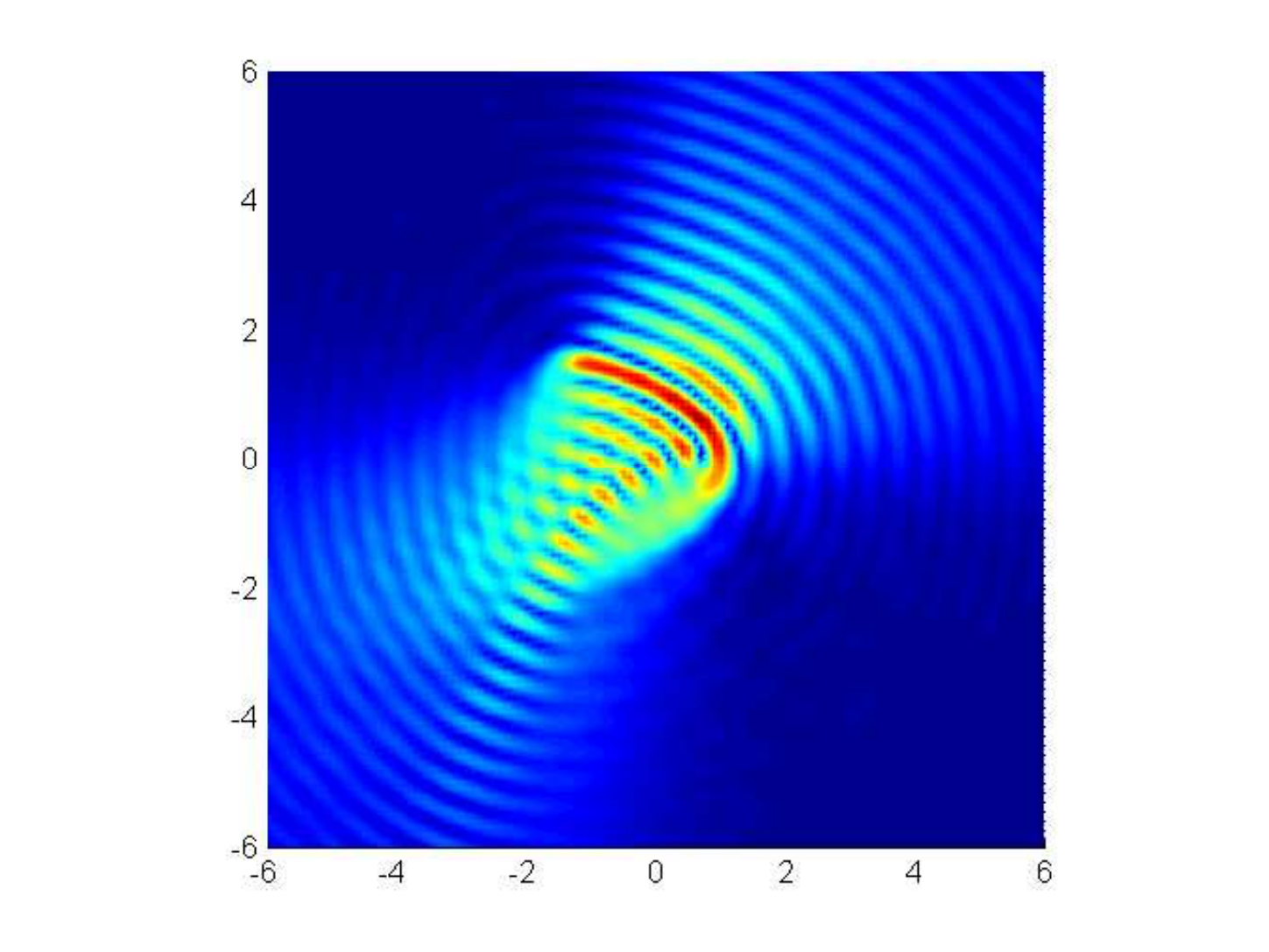}}\hspace*{-1cm}
  \subfigure[\textbf{$k=10$}]{
    \includegraphics[width=2in]{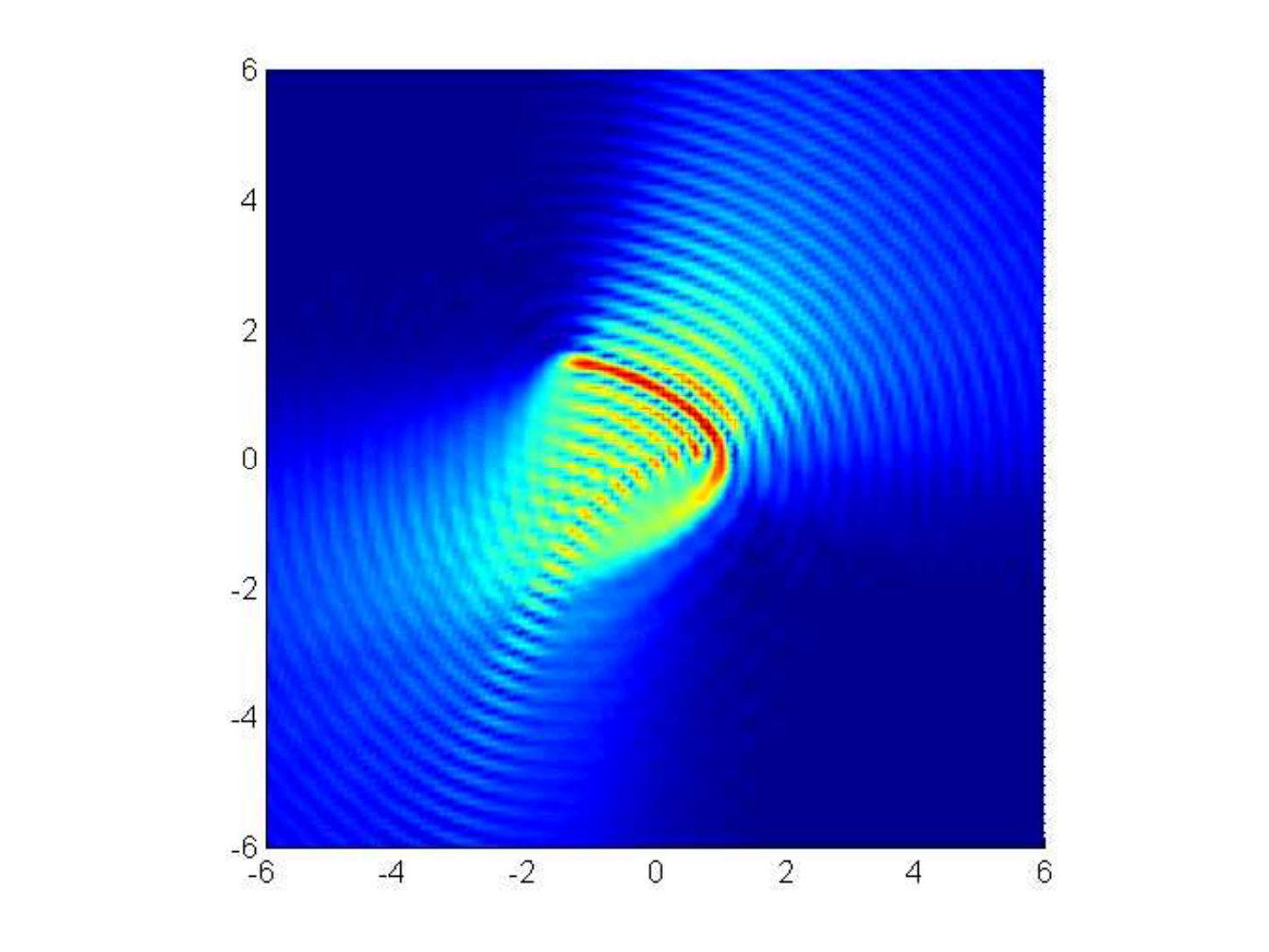}}
\caption{Reconstructions of the kite shaped obstacle with observation directions in $(0,\pi/2)$, $10\%$ noise and different wave numbers. }
\label{wavenumbers}
\end{figure}

Denote by $\hat{x}:=(\cos\phi,\sin\phi)$ the observation direction
with angle $\phi\in (0,2\pi)$. In the first example we consider the
set $\{u^{\infty}(\hx,\hth): \hx\in \Gamma, \hth\in
\mathbb{S}^{1}\}$ of far field patterns as data for the inverse
problem, i.e., a limited observation aperture $\Gamma\subset \mathbb{S}^{1}$
and the full incident aperture $\Sigma= \mathbb{S}^{1}$. In Figures
\ref{kiteapertures}-\ref{kiteapertures2}, we show the
reconstructions of a kite shaped obstacle with different limited
observation aperture $\Gamma$. As in many other numerical methods,
see e.g. \cite{L4}, a typical feature of the limited-aperture
problems is that the "shadow region" is elongated in down range.
However, Figures \ref{kiteapertures}-\ref{kiteapertures2} show our
results and confirm that, despite the use of the limited aperture data,
the region $\pa D_{+}(\hx):=\{y\in\pa D| \,\nu(y)\cdot\hat{x}>0\}$
with $\nu\in\mathbb{S}^1$ signifying the exterior unit normal vector
to $\partial D$, which can be observed directly from the direction
$\hat{x}$, is quite well captured. The imaging functional
\eqref{IndicatorDSM} produces a localized reconstruction of the
obstacle. As seen in Figure \ref{kiteapertures2}, the quality of the
reconstruction deteriorates as the aperture decreases. Physically,
the information from the "shadow region" $\pa
D_{-}(\hat{x}):=\{y\in\pa D |\, \nu(y)\cdot \hat{x}<0\}$ is very
weak, especially for high frequency waves \cite{MagerBleistein}. We
also refer to Figure \ref{wavenumbers} for the reconstructions with
the same observation aperture, but different wave numbers.

In the following, we provide a partial theoretical explanation to
the localized features of the reconstructions in the above numerical
examples. For a sound-soft obstacle $D$, the far field pattern of
the scattered field $u^s$ is given by \ben u^{\infty}(\hat{x},\hat
\theta)=-\int_{\pa D}\frac{\pa u}{\pa\nu}(y)e^{-ik\hat{x}\cdot
y}ds(y),\quad \hat{x}\in \mathbb{S}^{n-1}. \enn For large wave numbers $k$,
i.e., where the Kirchhoff approximation is considered, the far-field
pattern is given by \cite{CK} \be\label{Kirchhoff1}
u^{\infty}(\hat{x},\hat\theta)=-2\int_{\pa D_{-}(d)}\frac{\pa
e^{ik\hat\theta\cdot y}}{\pa\nu(y)}e^{-ik\hat{x}\cdot y}ds(y),\quad
\hat{x}\in \mathbb{S}^{n-1}, \en where $\pa D_{-}(\hat\theta):=\{y\in\pa D
|\, \nu(y)\cdot \theta<0\}$ is the region illuminated by the plane
wave with the incident direction $\hat\theta$. Clearly, for a fixed
$\hat\theta\in \mathbb{S}^{n-1}$, if the illuminated part $\pa
D_{-}(\hat\theta)$ is known in advance or obtained approximately,
one may have the far-field data $u^{\infty}(\hat{x},\hat\theta)$ for
all $\hat{x}\in \mathbb{S}^{n-1}$ by the formula \eqref{Kirchhoff1}. Using
further the well-known reciprocity relation for the far-field
pattern, we have \be\label{Kirchhoff2}
u^{\infty}(\hat{x},\hat\theta)
&=&u^{\infty}(-\hat\theta,-\hat{x})\cr &=&-2\int_{\pa
D_{+}(\hat{x})}\frac{\pa e^{-ik\hat{x}\cdot
y}}{\pa\nu(y)}e^{ik\hat\theta\cdot y}ds(y),\quad \hat{x}\in \mathbb{S}^{n-1}.
\en The representation \eqref{Kirchhoff2} implies that the shadow
domain $\pa D_{-}(\hat{x})$ makes no contribution to the far-field
data $u^{\infty}(\hat{x},\hat\theta)$. This fact partially explains,
in Figures \ref{kiteapertures}-\ref{wavenumbers},  why only the
region $\pa D_{+}(\hat{x})$ can be well reconstructed.

\begin{figure}[!htbp]
  \centering
  \subfigure[\textbf{Backscattering data for $\phi\in(0,\pi)$}]{
    \includegraphics[width=2.3in]{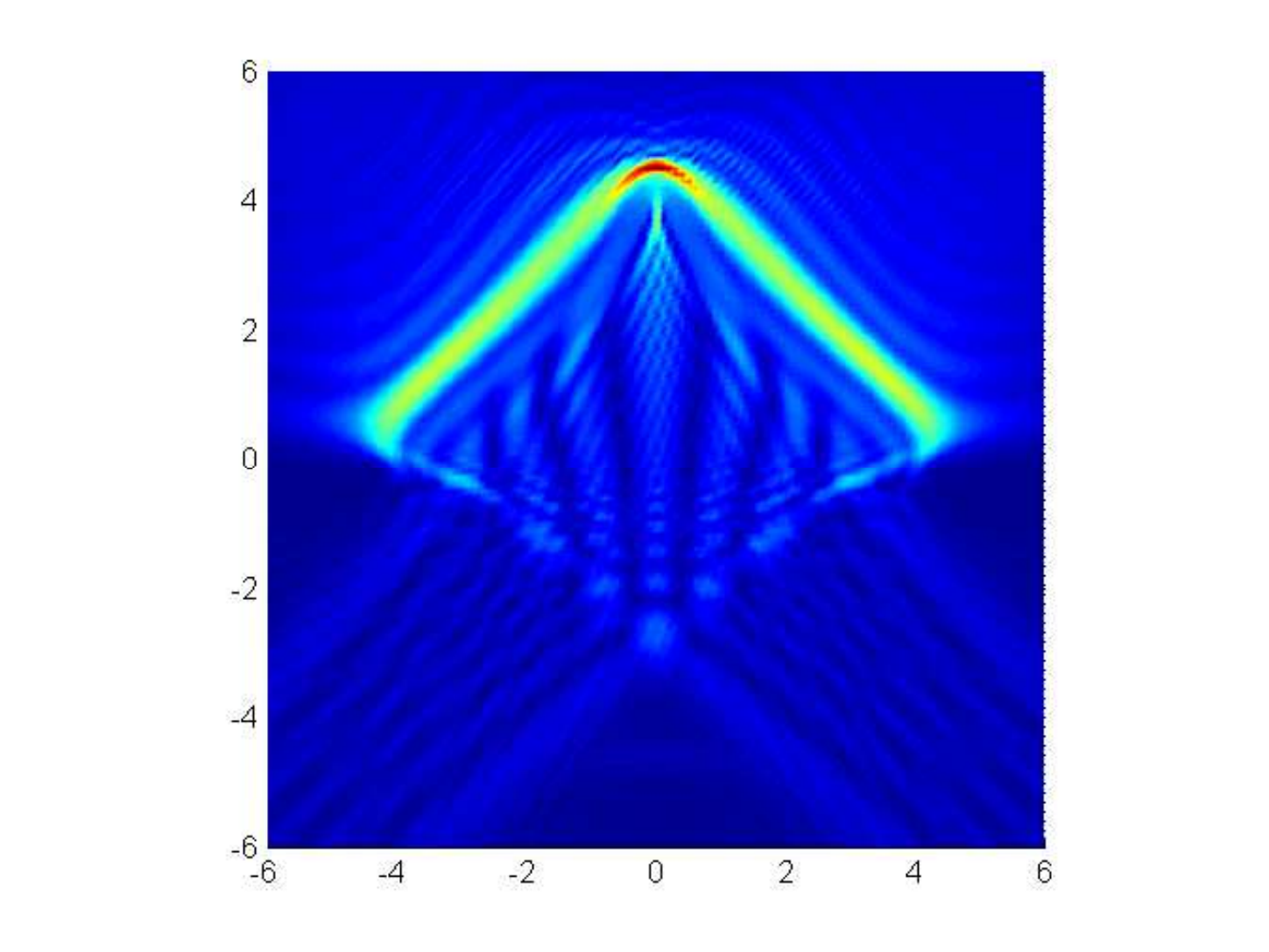}}
  \subfigure[\textbf{Backscattering data for $\phi\in(0,\pi/2)$}]{
    \includegraphics[width=2.3in]{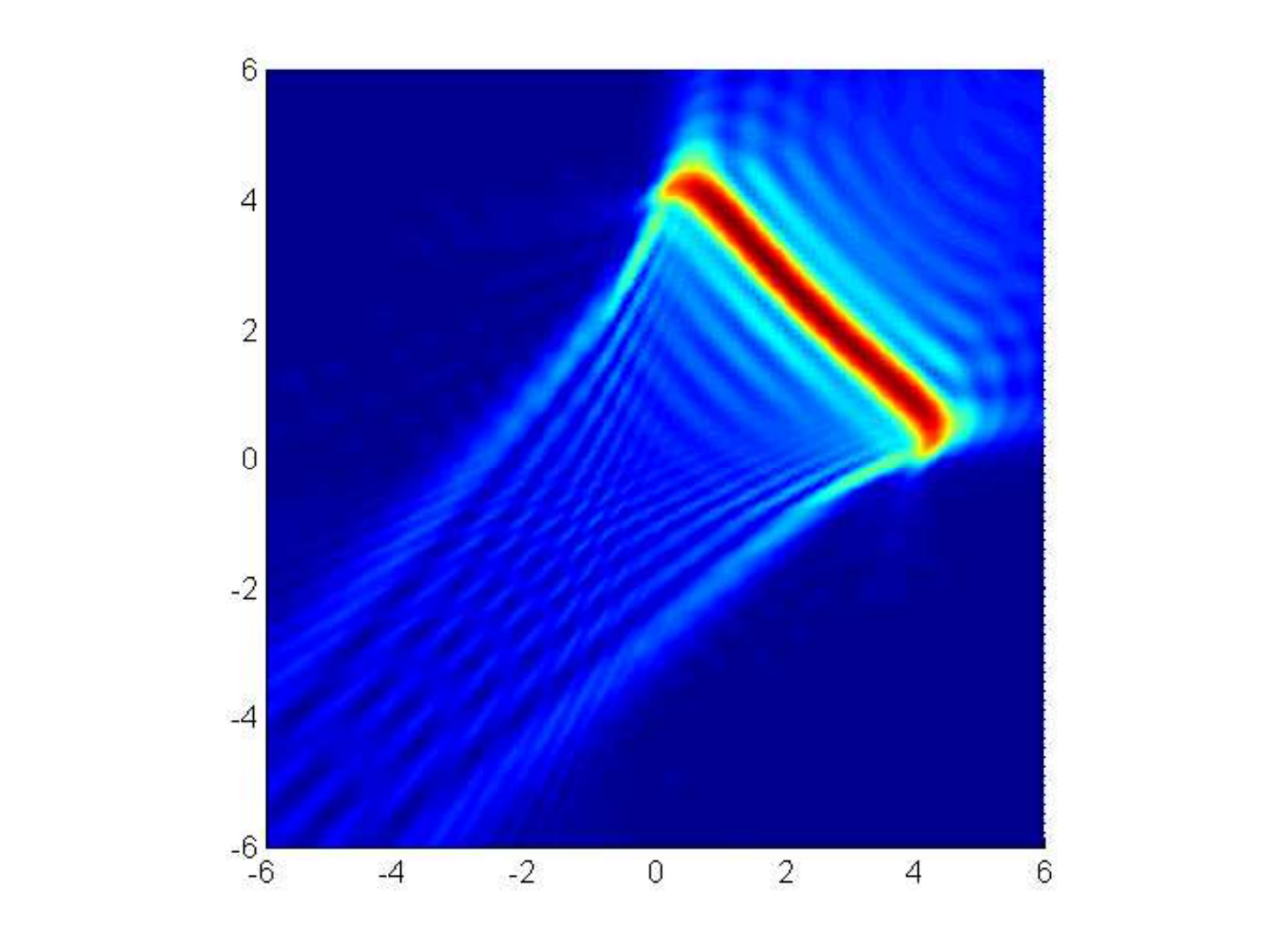}}
\caption{Reconstructions of a sound-soft round square by direct sampling method with limited aperture backscattering data, $k=20$, $10\%$ noise is added. }
\label{RoundSquare}
\end{figure}

We further consider the set $\{u^{\infty}(\hx,\hth): \hx\in \Gamma,
\hth\in \Sigma\}$ of far-field patterns as data for the inverse
problem \eqref{eq:ip1}, i.e., a limited observation aperture
$\Gamma\subset \mathbb{S}^{1}$ and a limited incident aperture $\Sigma\subset
\mathbb{S}^{1}$. Of particular interest is the case, where $\Gamma=-\Sigma$,
i.e., "backscattering" data is considered. Figure \ref{RoundSquare}
shows the corresponding reconstructions for a sound-soft round
square. It is clear that, similar to the previous examples, only the region $\pa
D_{+}(\hx)$ can be well captured and the quality of the
reconstruction deteriorates as the aperture decreases.
The numerical examples above show the intrinsic ill-posedness of
the inverse problem \eqref{eq:ip1} with limited aperture data due to
the lack of information. It is remarked that due to the localized
feature of the reconstruction method \eqref{IndicatorDSM}, if full
aperture data is available, say e.g., in the numerical examples in
Figure \ref{RoundSquare}, one can clearly recover the full obstacle
in an accurate manner. We would also like to emphasize that the
imaging functional \eqref{IndicatorDSM} does not work for the phaseless
inverse problem \eqref{eq:phaseless1}. That is, if only phaseless
measurement data is available, even with full aperture, one will
encounter even more severe ill-posedness due to the intrinsic lack
of information. Clearly, this intrinsic lack of information cannot
be remedied by any mathematical tricks unless more a-priori
information is complemented. This naturally leads us to integrate
the machine learning techniques into the inverse scattering problems
associated with incomplete measurement data.
In the next section, we shall briefly go through the main ingredients of the CNN for our
subsequent use.

%%%%%%%%%%%%%%%%%%%%%%%%%%%%%%%%%%%%%%%%%%%%%%%%%%%%%%%%%%%%%%%%%%

\section{Artificial neural network}

The concept of artificial neural network (ANN) was proposed by
McCulloch and Pitts in 1943, and they also proposed network
structure and mathematical descriptions of neurons through the M-P
model. In 1986, Rumelhart et. al. proposed the back propagation (BP)
algorithm for multi-layer perceptron and used the sigmoid function
as the nonlinear mapping \cite{RHW86}. This method solved the nonlinear
classification problem effectively, which led to a second
wave of neural networks development. In 1988, LeCun et al. began to
study convolutional neural networks (CNN) and developed well-known
LeNet5 \cite{LBDHHHJ90}, which is a special ANN with convolutional
kernel. It is pointed out that the BP algorithm has a gradient
disappearance problem which results in very poor learning efficiency.
In 2006, Hinton proposed a solution to the gradient
disappearance problem in deep network training, that is,
unsupervised pre-training initializes weights plus supervised
training fine-tuning \cite{HOT06,HS06}, which led to a third
wave of neural networks development. A deep neural network (DNN) is
an ANN with multiple layers between the input and output layers,
which is based on a hierarchy of composition of linear functions and
a given nonlinear activation function. It became popular and widely
accepted around 2010 due to the development of efficient learning
algorithms \cite{HOT06,HS06,KSH12,LBBH98,LKF10} and hardware
speed-ups such as the use of GPUs. In 2017, AlphaGo defeated Go world
champion Ke Jie by using Monte Carlo tree search combined with
two deep neural networks. These developments and applications show
that artificial neural networks have become very powerful tools in machine
learning or artificial intelligence such as image processing and
natural language processing \cite{GBC16,Zhou16}.

Mathematical analysis of neural networks has been developed by
several researchers. In 1989, Cybenko and Hornik have shown the
``universal approximation theorem'' independently, that is, under
mild assumptions on the activation function, a feed-forward network
with a single hidden layer can approximate continuous functions on
compact subsets of $\mathbb{R}^n$ \cite{Cyb89,HSW89}. However, this
theorem does not give the relationship between the convergence and the number of hidden layers and units,
and we refer to some extensions of the universal approximation
property \cite{GP90,GJP95,LLPS93}. On the other hand, Barron has
analyzed the two-layer neural networks in the Barron spaces around
1993, which is an extension of the Fourier analysis of two-layer
sigmoidal neural networks \cite{Bar93,Bar94}. In \cite{SCC18},
Shaham et. al. constructed a wavelet frame for approximating
functions on smooth manifolds, where the wavelets are computed from
Rectified Linear Units (ReLU). Under this construction, they have
shown how the specified size depends on the complexity of the
function and the curvature of the manifold, which can specify the
network architecture to obtain desired approximation properties.
Based on \cite{Bar93,Bar94}, E, Ma, and Wu have establish direct and
inverse approximation theorems, for functions in the Barron space
for two-layer networks, and also for functions in compositional
function space for residual neural network models \cite{EMW19}.

With the fast development of optimization algorithms and hardware
configuration, the zoo of neural network types grows exponentially,
which includes CNN, recurrent neural network (RNN,\cite{RHW86}),
long short-term memory (LSTM,\cite{HS97}), DNN, and etc. Among
these different types of neural networks, CNNs are used heavily in
image processing. Firstly, the CNNs provide an advantage over
feed-forward networks because they are capable of considering
locality of features and extracting spatial features. The
convolution layer will transform an input into a stack of feature
mappings. The depth of the feature map stack depends on how many
filters are defined for a layer. Secondly, with the fixed small
window size of convolutions, a CNN minimizes the computation
sharply compared to a regular neural network.

Next, we present the concepts and notations of the CNN. A CNN
consists of $L+1$ layers, where the layer $0$ is the input layer,
the layers $0<l<L$ are the hidden layers, and the layer $L$ is the output
layer. We can adopt the $\mathbb{F}_{12}$ and $\mathbb{F}_{full}$ defined in what follows in (\ref{eq:Full})
as input and output layer respectively to retrieval the full-aperture far-field directly.
The activation functions in the hidden layers can be any
activation function such as sigmoids, rectified linear units, or
hyperbolic tangents. Here, we will use parametric rectified linear
unit (PReLU) in the hidden layers.

Parametric Rectified Linear Unit (PReLU) is a type of leaky
Rectified Linear Unit (ReLU),  where instead of having a
predetermined slope such as 0.01. The slope is a parameter to be
determined by the neural network
\begin{equation*}
    \sigma(\alpha, x)=
    \begin{cases}
        \alpha x, & \; x < 0, \\
        x, & \; x \geq 0.
    \end{cases}
\end{equation*}
where $\alpha$ is a learned array with the same dimension as $x$.

Convolution layer is composed of several convolution kernels which
are used to compute different feature maps. The complete feature
maps are obtained by using several different kernels. Let ${\bf
z}^{l}_{i,j}$ be the input patch centered at location $(i, j)$ of
the $l$-th layer and $y^{l}_{i,j,k}$ be the feature value at
location $(i, j)$ in the $k$-th feature map of the $l$-th layer,
respectively. Then
\begin{equation*}
y^{l}_{i,j,k}= ({\bf w}^{l}_{k})^{T} {\bf z}^{l}_{i,j}+b^{l}_k,
\end{equation*}
where ${\bf w}^{l}_{k}\in \,\mathbb{R}^{M\times M}$ and $b^{l}_k$ are the weight vector with window size $M\times M$ and bias
term of the $k$-th filter of the $l$-th layer, respectively. Note
that the kernel ${\bf w}^{l}_{k}$, which generates the feature map
$y^{l}_{:,:,k}$, is shared. Such a weight sharing mechanism has
several advantages. For example, it can reduce the model complexity
and make the network easier to train. The activation function in
layer $l$ will be denoted by $\sigma^{l}=\sigma(\alpha, x)$. Then
the activation value $z^{l}_{i,j,k}$ of convolutional feature
$y^{l}_{i,j,k}$ can be computed as
\begin{equation*}
z^{l}_{i,j,k}=\sigma^{l}(y^{l}_{i,j,k}).
\end{equation*}
Denote ${\bf w}=\{{\bf w}^{l}_{k}\}$ and $b=\{b^{l}_k \}$. The
schematic representation of general convolution neural network is
shown in Fig~\ref{fig:schematic representation CNN}.
\begin{figure}[!htbp]
\centering
\includegraphics[width=.7\textwidth]{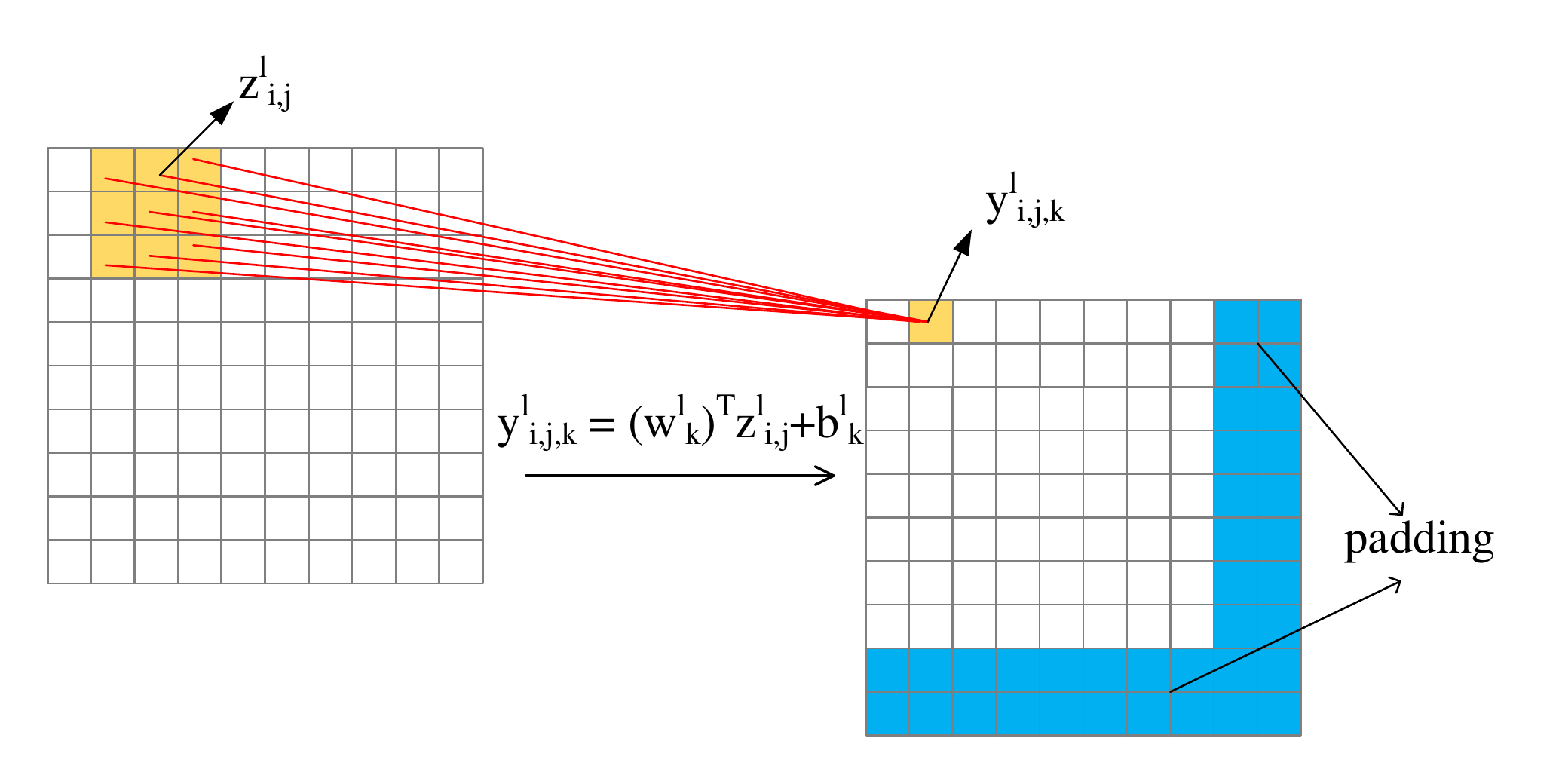}
\caption{Schematic representation of general
convolution neural network.} \label{fig:schematic representation
CNN}
\end{figure}

%\newpage

%%%%%%%%%%%%%%%%%%%%%%%%%%%%%%%%%%%%%%%%%%%%%%%%%%%%%%%%%%%%%%%%%%%%%%%%

\section{Algorithmic development and numerical experiments}

%In this section, we present three numerical simulations to verify the efficiency of our
%proposed method, i.e., convolution neural networks (CNN) for the
%recovery of the full data matrix from limited measure data. We shall
%introduce some preliminary concepts for the sequel numerical
%simulations.

In this section, we present our data retrieval method based on CNN for the inverse scattering problems \eqref{eq:ip1} and \eqref{eq:phaseless1}. The proposed method shall first try to recover the full data by using the provided incomplete data via certain training mechanisms through proper CNNs. After the recovery of the full data, we then make use of the imaging method \eqref{IndicatorDSM} to reconstruct the unknown obstacle. We would like to remark that this data retrieval and inverse scattering recovery process can obviously be extended to the other inverse scattering problems with incomplete data. Moreover, mainly due to limiting computing resources, we only consider 2D examples in what follows, and the extension to three dimensions are obvious if more computing resources are available.

We introduce some general ingredients for our
subsequent algorithmic developments through concrete examples.
Denote by $\mathcal{X}$ the set of far-field patterns on
$\Gamma \times \Sigma$ and $\mathcal{Y}$ the set of far-field patterns
on $\mathbb{S}^{1} \times \mathbb{S}^{1}$.
Then each limited-aperture far-field pattern $X \in \mathcal{X}$ is
mapped to the unique full-aperture far-field pattern
\[
  Y=\mathscr{A}(X) \in \mathcal{Y},
\]
where $\mathscr{A}$ denotes the analytic continuation.
Our objective is to retrieve the function $\mathscr{A}: \mathcal{X} \to \mathcal{Y}$ from the samples
\begin{equation*}
  \mathcal{S} = \left\{ (X^n \in \mathcal{X},Y^n \in \mathcal{Y}) \right\}_{n=1}^{n^*},
\end{equation*}
where $X^n$, $Y^n$ are a limited-aperture and the full-aperture far-field data
corresponding to the same obstacle $\Omega_n$ from the sample set
\[
  \mathcal{D}=\{\Omega_n\}_{n=1}^{n^*},
\]
i.e.,
\begin{align*}
  X^n &= u^\infty(\Omega_n; \hat{x}, \hat{\theta}), \quad (\hat{x},\hat{\theta}) \in \Gamma \times \Sigma,\\
  Y^n &= u^\infty(\Omega_n; \hat{x}, \hat{\theta}), \quad (\hat{x},\hat{\theta}) \in \mathbb{S}^{1} \times \mathbb{S}^{1}.
\end{align*}
Here, $n^*$ denotes the total number of samples.

It is generally an open problem as how to formulate and generate a set of random and arbitrary non-selfintersecting shapes. In this paper, we restrict ourself to star
shaped domains whose boundaries are defined by the parametric curves
\begin{equation*}
  x(t) = r(t) \cos (t), \quad y(t) = r(t) \sin (t), \quad t \in [0,2\pi].
\end{equation*}
We then represent the radius function \( r(t) \) by the truncated Fourier series expansion
\begin{equation*}
  r(t) = a_0 \left\{ 1 + \frac{1}{2N} \sum_{n=1}^N n^{-q} \left[ a_n \cos (nt) + b_n \sin (nt) \right]
  \right\},
\end{equation*}
where \( N \in \mathbb{N} \) is the cut-off frequency, and $a_n,b_n, n=1,\cdots,N$ are random numbers draw from
the uniform distribution in $[-1,1]$. The factor \( \dfrac{ 1 }{ 2N } \) is to guarantee $r(t)$ does not change sign so that the resultant shape is non-selfintersecting. The parameter \( q \) can be adjusted to control the decreasing rate of the Fourier coefficients and thus the level of smootheness of the shape. The parameter \( a_0 \) can be taken as a fixed value or a random variable.

The range of the sample set can be adjusted by the parameter \( N \) and the range of \( a_0 \). A broader range of the sample set implies a broader applicability of the training results, but also requires a larger training set and training time. In the subsequent numerical experiments, we take \( N=5 \), \( q=0 \), and \( a_0 \) to be a random number drawn from the uniform distribution in \( [0.5, 1.5] \).

We also need to introduce the multi-static response (MSR) matrix and
the related properties. The multi-static response matrix
$\mathbb{F}_{full} \in C^{2m\times 2m}$ corresponding to the full-aperture far-field data is defined as
\begin{eqnarray*} \mathbb{F}_{full}=\left(
\begin{array}{cccc}\vspace{2mm}
u_{1,1}^{\infty} &  u_{1,2}^{\infty} & \cdots & u_{1,2m}^{\infty} \\
u_{2,1}^{\infty} &  u_{2,2}^{\infty} & \cdots & u_{2,2m}^{\infty} \\
\vdots &  \vdots & \ddots & \vdots \\
u_{2m,1}^{\infty} &  u_{2m,2}^{\infty} & \cdots & u_{2m,2m}^{\infty}
\end{array}
\right),
\end{eqnarray*}
where $u_{i,j}^{\infty}=u^{\infty}(\hat{x}_j; \hat{\theta}^i)$ for
$1 \leq i, j \leq 2m$ corresponding to $2m$ observation directions
$\hat{x}_j$ and $2m$ incident directions $\hat{\theta}^i$. Note that
$\mathbb{F}_{full}$ is neither symmetric nor Hermitian, i.e.,
\begin{eqnarray*}
\mathbb{F}_{full}\neq \mathbb{F}_{full}^{T}\quad and \quad
\mathbb{F}_{full}\neq \mathbb{F}_{full}^{*}.
\end{eqnarray*}
Here and throughout the paper we use the superscript ``T'' and ``*''
to denote the transpose and the conjugate transpose, respectively,
of a matrix.

Generally speaking, we can partition the $2m$-by-$2m$ MSR matrix
$\mathbb{F}_{full}$ into a 2-by-2 block matrix
\begin{eqnarray}\label{eq:Full} \mathbb{F}_{full}=\left(
\begin{array}{cccc}\vspace{2mm}
\mathbb{F}_{11} & \mathbb{F}_{12}  \\
\mathbb{F}_{21} & \mathbb{F}_{22}
\end{array} \right),
\end{eqnarray}
where $\mathbb{F}_{11}\in C^{m_1\times m_1}$, $\mathbb{F}_{12}\in
C^{m_1\times m_2}$, $\mathbb{F}_{11}\in C^{m_2\times m_1}$,
$\mathbb{F}_{12}\in C^{m_2\times m_2}$, and $m_1+m_2=2m$. Here,
$\mathbb{F}_{12}$ denotes the limited-aperture far-field data.

Let $n_1$ and $n_2$ be the size of train set and test set,
respectively. In the sequel, we denote $\mathbb{F}^{n}_{full}$ by
the $n$-th sample in the sample set, and similar notations for sub-block
matrices $\mathbb{F}^{n}_{ij}$, $i,j=1,2$.

We proceed to the algorithmic development through four concrete examples.

\vspace{2mm}

%%%%%%%%%%%%%%%%%%%%%%%%%%%%%%%%%%%%%%%%%%%%%%%%%%%%%%%%%%%%%%%%%%%%%%%

%%%%%%%%%%%%%%%%%%%%%%%%%%%%%%%%%%%%%%%%%%%%%%%%%%%%%%%%%%%%%%%%%%%%%%%

\noindent {\bf Example 1.} In this example, we consider the standard
loss function with the total number of samples $n^{*}=20000$.

\begin{figure}[!htbp]
\centering
\includegraphics[width=\textwidth]{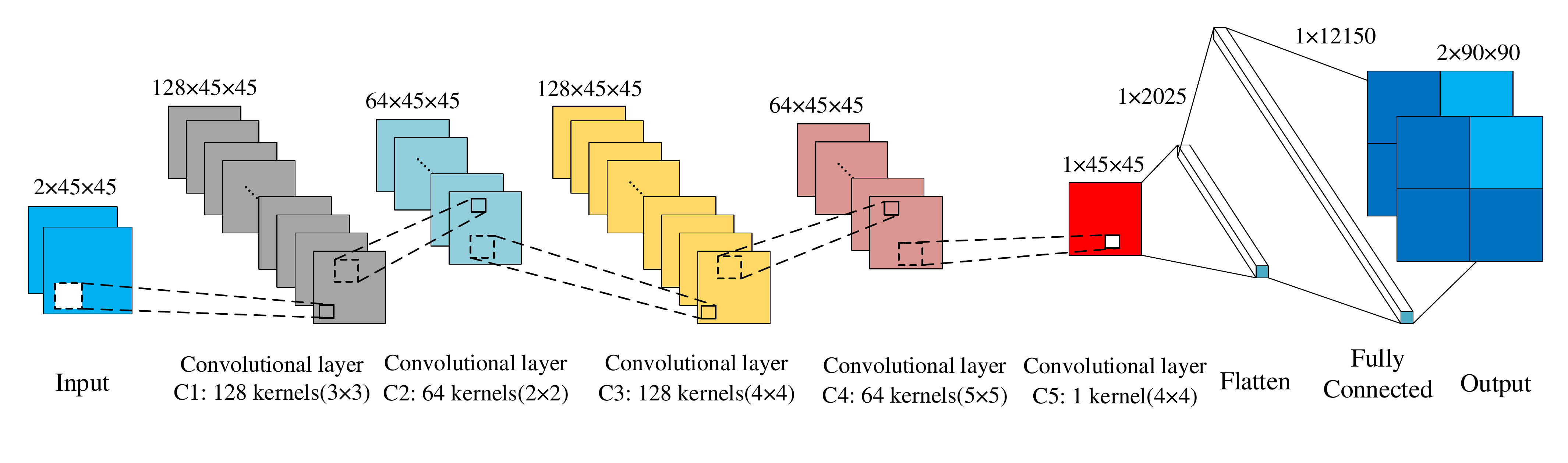}
\caption{The architecture of our convolution neural network in
Example 1 and 2.} \label{fig:EXAM12CNN}
\end{figure}

Following \cite{SD19}, the architecture of our neural network is a
feed-forward stack of five sequential combinations of the
convolution, batch normalization and PReLU layers, followed by one
fully connected layer. The numbers of filters in the five
convolutional layers are 128, 64, 128, 64 and 1, respectively, which
are shown in Fig~\ref{fig:EXAM12CNN}. The corresponding window sizes
of convolutions are $3\times 3$, $2\times 2$, $4\times 4$, $5\times
5$, and $4\times 4$, respectively(i.e. $M=3,2,4,5$, and $4$, respectively for each layer).
The stride of the convolution is one, and zero-padding is used. The initial value of the bias is
zero. The weight initialization is via the Glorot uniform
initializer \cite{GB10}.

With this architecture, we train the network with the Adam
optimizer, using a mini-batch of 128 samples in each iteration, and
200 epochs. The initial learning rate and forgetting rate of the
Adam are the same as in \cite{KB14}. The total number of samples is
$n^{*}=20000$, where the size of the train set is $n_1=16000$ and the
size of the test set is $n_2=4000$. For each sample in the training set
$(n=1,\cdots,n_1)$, the input data is $\mathbb{F}^{n}_{12}\in
C^{m_1\times m_1}$ (with $m_1=45$ in first three simulations) which
is the upper-right quarter of $\mathbb{F}^{n}_{full}$, while the target
data are $\mathbb{F}^{n}_{11}$, $\mathbb{F}^{n}_{12}$,
$\mathbb{F}^{n}_{21}$, and $\mathbb{F}^{n}_{22}$. That is, our CNN
defines a mapping $ \mathbb{R}^{m_1\times m_1 \times 2}\rightarrow
\mathbb{R}^{2m\times 2m \times 2}$ (each complex number has a real
part and an imaginary part).

Given the input data $x=\mathbb{F}^{n}_{12}\in \mathbb{R}^{m_1\times
m_1\times 2}$ and the target outputs
$y=[\mathbb{F}^{n}_{11},\mathbb{F}^{n}_{12};\mathbb{F}^{n}_{21},\mathbb{F}^{n}_{22}]\in
\mathbb{R}^{2m\times 2m \times 2}$ in the training set, we wish to
choose our weights and biases such that $y^{C}(x;{\bf w},b)$ is a
good approximation of $y(x)$, where $y^{C}(x;{\bf
w},b)=[\mathbb{F}^{n,C}_{11},\mathbb{F}^{n,C}_{12};\mathbb{F}^{n,C}_{21},\mathbb{F}^{n,C}_{22}]$
denotes the outputs of CNN.

To find the weights and biases, we define the following loss
function
\begin{eqnarray}\label{eq:L1 def}
\mathcal{L}_1&=&\mathcal{L}_1(y,y^{C})
=\frac{1}{n_1}\sum\limits_{n=1}^{n_1}\sum\limits_{i,j=1}^{2}
\|\mathbb{F}^{n}_{ij}-\mathbb{F}_{ij}^{n,C}\|^2_{F},
\end{eqnarray}
where $\|\cdot\|_{F}$ is the Frobenius norm of a matrix. We use the
CNN to compute
\begin{eqnarray*}
{\bf w}^*, b^*&=& \mbox{arg}\min\limits_{{\bf w},b} \mathcal{L}_1(y,
y^{C})
\end{eqnarray*}
with the Adam optimizer. We define the following sub-block relative
errors
\begin{eqnarray*}
e^{n}_{11}=\frac{\|\mathbb{F}^{n}_{11}-\mathbb{F}_{11}^{n,C}\|_{F}}{\|\mathbb{F}^{n}_{11}\|_{F}},~
e^{n}_{21}=\frac{\|\mathbb{F}^{n}_{21}-\mathbb{F}_{21}^{n,C}\|_{F}}{\|\mathbb{F}^{n}_{21}\|_{F}},~
e^{n}_{22}=\frac{\|\mathbb{F}^{n}_{22}-\mathbb{F}_{22}^{n,C}\|_{F}}{\|\mathbb{F}^{n}_{22}\|_{F}},
\quad n=n_1+1,\cdots,n^*
\end{eqnarray*}
for each sample in the test set, and use the sub-block relative errors
of the test set
\begin{eqnarray*}
\bar{e}_{11}=\frac{1}{n_2}\sum\limits_{n=n_1+1}^{n^*}e^{n}_{11},~~
\bar{e}_{21}=\frac{1}{n_2}\sum\limits_{n=n_1+1}^{n^*}e^{n}_{21},~~
\bar{e}_{22}=\frac{1}{n_2}\sum\limits_{n=n_1+1}^{n^*}e^{n}_{22},
\end{eqnarray*}
to check the efficiency of our CNN. Table~\ref{table:sub_block_ex1} shows that the sub-block
relative errors of the test set are decreasing as $n_1$ and $n_2$
increase. We emphasize that the errors $\bar{e}_{11}$ and
$\bar{e}_{22}$ are smaller than $\bar{e}_{21}$, which is reasonable
since $\mathbb{F}_{11}$ and $\mathbb{F}_{22}$ are closer to the
$\mathbb{F}_{12}$.

\begin{table}[!htbp]
  \caption{The sub-block relative errors of test set in Example 1.}\label{table:sub_block_ex1}
\begin{tabular}{|r||c|c|c|c|c|c|c|c|c|}\hline
              & \multicolumn{3}{|c|}{error real(\%) } &   \multicolumn{3}{|c|}{error image(\%)} &
    \multicolumn{3}{|c|}{error norm(\%)}\\ \hline
   $(n_1,n_2)$   &  $\bar{e}_{11}$   &  $\bar{e}_{21}$   &  $\bar{e}_{22}$       &  $\bar{e}_{11}$& $\bar{e}_{21}$ & $\bar{e}_{22}$ &$\bar{e}_{11}$ &$\bar{e}_{21}$ & $\bar{e}_{22}$\\ \hline \hline
   (2000,500)    &   26.01   &    84.15       &  26.01   & 30.72   & 86.03  & 30.71   & 28.06   & 84.99   & 30.71 \\ \hline
   (4000,1000)   &   19.01   &    67.55       &  19.00   & 22.49   & 68.43  & 22.49   & 20.52   & 67.92   & 20.52 \\ \hline
   (8000,2000)   &   14.22   &    51.15       &  14.22   & 16.88   & 52.12  & 16.88   & 15.38   & 51.54   & 15.38\\ \hline
   (16000,4000)  &   10.66   &    38.88       &  10.66   & 12.61   & 39.82  & 12.61   & 11.50   & 39.27   & 11.50 \\ \hline
\end{tabular}
\end{table}

For $n=n_1+1,\cdots,n^*$, combining the outputs
$[\mathbb{F}^{n,C}_{11},\mathbb{F}^{n,C}_{21},\mathbb{F}^{n,C}_{22}]$
and the input $\mathbb{F}^{n}_{12}$ of CNN, we can obtain the
recovery multi-static response matrix as
\begin{eqnarray*}
  \mathbb{\widetilde{F}}^{n}_{full}=\left(
\begin{array}{cccc}\vspace{2mm}
\mathbb{F}^{n,C}_{11} & \mathbb{F}^{n}_{12}  \\
\mathbb{F}^{n,C}_{21} & \mathbb{F}^{n,C}_{22}
\end{array} \right), \quad n=n_1+1,\cdots,n^*,
\end{eqnarray*}
for each sample in the test set. Given a noise-free $2m\times 2m$ sample
$\mathbb{F}^{n}_{full}$ in the test set and its approximation
$\mathbb{\widetilde{F}}^{n}_{full}$, define
\begin{eqnarray*}
e^{n}&=&\frac{\|\mathbb{F}^{n}_{11}-\mathbb{F}_{11}^{n,C}\|_{F}+\|\mathbb{F}^{n}_{21}-\mathbb{F}_{21}^{n,C}\|_{F}+\|\mathbb{F}^{n}_{22}-\mathbb{F}_{22}^{n,C}\|_{F}}
{\|\mathbb{F}^{n}_{11}\|_{F}+\|\mathbb{F}^{n}_{21}\|_{F}+\|\mathbb{F}^{n}_{22}\|_{F}},\\
\mbox{MSE}^{n}&=&{\frac{1}{2m\cdot
2m}}\sum_{i=1}^{2m}\sum_{j=1}^{2m}
\left|\mathbb{F}^{n}_{full}(i,j)-\mathbb{\widetilde{F}}^{n}_{full}(i,j)\right|^{2},\\
\mbox{PSNR}^{n}&=&10\cdot \mbox{log}_{10}
\left(\frac{\max\limits_{i,j}|\mathbb{F}^{n}_{full}(i,j)|}{\mbox{MSE}^{n}}\right).
\end{eqnarray*}
In order to check the efficiency of our CNN, we also introduce the
relative error, mean square error (MSE), and peak signal-to-noise
ratio (PSNR) of the test set as follows:
\begin{eqnarray*}
\overline{e}=\frac{1}{n_2}\sum\limits_{n=n_1+1}^{n^*}e^{n},~
\overline{\mbox{MSE}}=\frac{1}{n_2}\sum\limits_{n=n_1+1}^{n^*}\mbox{MSE}^{n},~
\overline{\mbox{PSNR}}=\frac{1}{n_2}\sum\limits_{n=n_1+1}^{n^*}\mbox{PSNR}^{n}.
\end{eqnarray*}
Table~\ref{table:ex1} shows that the relative errors and $\overline{\mbox{MSE}}$s
of the test set are decreasing as $n_1$ and $n_2$ increase. Meanwhile,
the $\overline{\mbox{PSNR}}$s of test set are increasing as $n_1$
and $n_2$ increase, which means the more samples in the train test,
the better the prediction will be.
\begin{table}[!htbp]
  \caption{The relative errors, $\overline{\mbox{MSE}}$s, and $\overline{\mbox{PSNR}}$s of test set in Example 1.}\label{table:ex1}
\begin{tabular}{|r||c|c|c|c|c|c|c|c|c|}\hline
             & \multicolumn{3}{|c|}{relative error $\overline{e}$(\%) } &   \multicolumn{3}{|c|}{$\overline{\mbox{MSE}}$($10^{-2}$)} &
    \multicolumn{3}{|c|}{$\overline{\mbox{PSNR}}$}\\ \hline
   $(n_1,n_2)$   &    real   &    image       &  norm    &  real   & image  & norm    & real    &image    & norm \\ \hline \hline
   (2000,500)    &   39.13   &    44.03       &  41.40   & 9.372   & 9.384  & 5.430   & 19.21   & 16.22   & 22.37 \\ \hline
   (4000,1000)   &   34.06   &    38.47       &  36.11   & 6.979   & 7.042  & 4.371   & 20.28   & 17.27   & 22.24 \\ \hline
   (8000,2000)   &   23.97   &    27.09       &  25.41   & 3.445   & 3.440  & 2.594   & 23.24   & 20.22   & 25.65\\ \hline
   (16000,4000)  &   16.63   &    18.80       &  17.63   & 1.691   & 1.692  & 1.365   & 26.46   & 23.46   & 28.63 \\ \hline
\end{tabular}
\end{table}

\begin{figure}[!htbp]
\centering
\includegraphics[width=\textwidth]{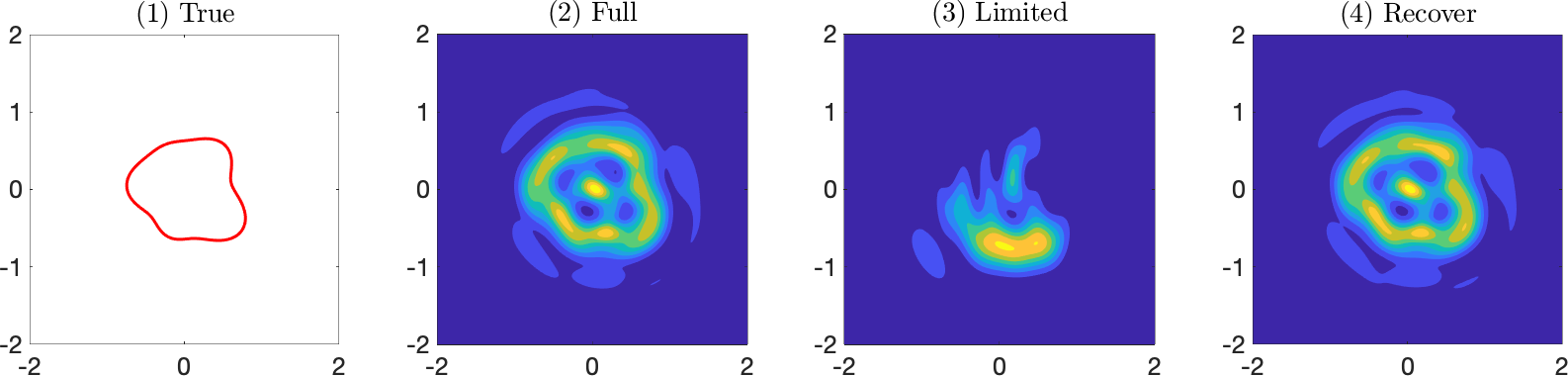}
\caption{Numerical constructions of the 1st shape in the test set. (1) the true shape; (2) reconstructed shape from the full MSR matrix
$\mathbb{F}^{n}_{full}$; (3) reconstructed shape from the limited MSR matrix
$\mathbb{F}^{n}_{12}$; (4) reconstructed shape from the recovered MSR matrix
$\mathbb{\widetilde{F}}^{n}_{full}$ via CNN.}
\label{fig:pred_rand}
\end{figure}

Fig~\ref{fig:pred_rand} shows the numerical constructions of the 1st shape from the test set by the full MSR matrix $\mathbb{F}^{n}_{full}$, a limited MSR
matrix $\mathbb{F}^{n}_{12}$, and the recovery MSR matrix
$\mathbb{\widetilde{F}}^{n}_{full}$ via CNN. Although the average
error $\bar{e}_{21}$ is around $20\%$, the numerical construction by the
recovery MSR matrix via CNN is indistinguishable from the numerical
construction by the full MSR matrix.

\vspace{3mm}

%\newpage

%%%%%%%%%%%%%%%%%%%%%%%%%%%%%%%

\noindent {\bf Example 2.} In this example, we consider the loss
function with a regularization term. The architecture of convolution
neural network, the input data $x=\mathbb{F}^{n}_{12}\in
\mathbb{R}^{m_1\times m_1 \times 2}$, and the target outputs
$y=[\mathbb{F}^{n}_{11},\mathbb{F}^{n}_{12};\mathbb{F}^{n}_{21},\mathbb{F}^{n}_{22}]\in
\mathbb{R}^{m\times m \times 2}$ are the same as in Example 1. The total
number of samples is $n^{*}=50000$, where the size of the training set is
$n_1=40000$ and the size of the test set is $n_2=10000$.

Define the modified loss function as
\begin{eqnarray*}
&&\mathcal{L}_2(y,y^{C})=\mathcal{L}_1(y,y^{C}) \\
&+&\frac{\alpha}{n_1}\sum\limits_{n=1}^{n_1}
\left\{|\mathbb{F}_{11}^{n,C}\cap
\mathbb{F}_{12}^{n}|^2_{H^1}+|\mathbb{F}_{21}^{n,C}\cap
\mathbb{F}_{22}^{n,C}|^2_{H^1} +|\mathbb{F}_{11}^{n,C}\cap
\mathbb{F}_{21}^{n}|^2_{H^1}+|\mathbb{F}_{12}^{n}\cap
\mathbb{F}_{22}^{n}|^2_{H^1}\right\}, \nonumber
\end{eqnarray*}
where $\mathcal{L}_1$ is specified in (\ref{eq:L1 def}), the second
summation is the regularization term on the interfaces of each
sub-block matrix, and the regularization parameter $\alpha=10^{-3}$.

Table~\ref{table:50T_ex2} and Table~\ref{table:BV_ex2} show that the relative errors and
$\overline{\mbox{MSE}}$s of the test set are decreasing as $n_1$ and
$n_2$ increase, while the $\overline{\mbox{PSNR}}$s of the test set
are increasing as $n_1$ and $n_2$ increase. Meanwhile the the
relative errors and $\overline{\mbox{MSE}}$s of $\mathcal{L}_2$ are
smaller than the $\mathcal{L}_1$, and $\overline{\mbox{PSNR}}$s of
$\mathcal{L}_2$ are great than the $\mathcal{L}_1$. These mean the
regularization term reduces the prediction error efficiently.
\begin{table}[!htbp]
  \caption{The relative errors, $\overline{\mbox{MSE}}$s, and $\overline{\mbox{PSNR}}$s of test set in Example 2 without regularization term.}\label{table:50T_ex2}
\begin{tabular}{|r||c|c|c|c|c|c|c|c|c|}\hline
             & \multicolumn{3}{|c|}{relative error $\overline{e}$(\%) } &   \multicolumn{3}{|c|}{$\overline{\mbox{MSE}}$($10^{-2}$)} &
    \multicolumn{3}{|c|}{$\overline{\mbox{PSNR}}$}\\ \hline
   $(n_1,n_2)$   &    real   &    image       &  norm    &  real   & image  & norm    & real    &image    & norm \\ \hline \hline
   (5000,1250)    &   29.91   &    33.58       &  31.62   & 5.489   & 5.511  & 3.650   & 21.51   & 18.59   & 24.31 \\ \hline
   (10000,2500)   &   21.54   &    24.32       &  22.83   & 2.878   & 2.894  & 2.130   & 24.33   & 21.31   & 26.62 \\ \hline
   (20000,5000)   &   15.67   &    17.74       &  16.63   & 1.438   & 1.439  & 1.198   & 26.85   & 23.84   & 28.95\\ \hline
   (40000,10000)  &   12.22   &    13.85       &  12.97   & 0.818   & 0.816  & 0.712   & 28.84   & 25.82   & 30.80 \\ \hline
\end{tabular}
\end{table}

\begin{table}[!htbp]
  \caption{The relative errors, $\overline{\mbox{MSE}}$s, and $\overline{\mbox{PSNR}}$s of test set in Example 2 with regularization term.}\label{table:BV_ex2}
\begin{tabular}{|r||c|c|c|c|c|c|c|c|c|}\hline
           & \multicolumn{3}{|c|}{relative error $\overline{e}$(\%) } &   \multicolumn{3}{|c|}{$\overline{\mbox{MSE}}$($10^{-2}$)} &
    \multicolumn{3}{|c|}{$\overline{\mbox{PSNR}}$}\\ \hline
   $(n_1,n_2)$    &    real   &    image       &  norm    &  real   & image  &  norm   & real    &image    & norm   \\ \hline \hline
   (5000,1250)    &   30.08   &    33.63       &  33.73   & 5.372   & 5.319  & 3.686   & 21.34   & 18.47   & 23.91  \\ \hline
   (10000,2500)   &   21.71   &    24.43       &  22.97   & 2.869   & 2.852  & 2.128   & 24.16   & 21.18   & 26.64  \\ \hline
   (20000,5000)   &   15.26   &    17.17       &  16.14   & 1.390   & 1.383  & 1.149   & 27.12   & 24.18   & 29.30  \\ \hline
   (40000,10000)  &   12.00   &    13.45       &  12.67   & 0.798   & 0.786  & 0.677   & 29.02   & 26.10   & 31.20  \\ \hline
\end{tabular}
\end{table}

As an example, let's consider the square shape with corners at \( (1.5, 0) \), \( (0, 1.5) \), \( (-1.5, 0) \) and \( (0, -1.5) \). This is a star shape with the radius function given by
\[
  r(t) =
  \begin{cases}
    1 / (\sin t + \cos t), & \quad t \in [0, \frac{ \pi }{ 2 }],\\
    1 / (\sin t - \cos t), & \quad  t \in [\frac{ \pi }{ 2 }, \pi],\\
    1 / (-\sin t - \cos t), & \quad  t \in [\pi, \frac{ 3\pi }{ 2 }],\\
    1 / (-\sin t + \cos t), & \quad  t \in [\frac{ 3\pi }{ 2 }, 2\pi].
  \end{cases}
\]
One can verify the Fourier coefficients of \( r(t) \) lies within the range of possible random numbers for the training and test sets, but it is not one of the \( 20000 \) actual samples used in this experiment.

Denote the full MSR matrix $\mathbb{F}^{s}_{full}$
by \begin{eqnarray*} \mathbb{F}^{s}_{full}=\left(
\begin{array}{cccc}\vspace{2mm}
\mathbb{F}^{s}_{11} & \mathbb{F}^{s}_{12}  \\
\mathbb{F}^{s}_{21} & \mathbb{F}^{s}_{22}
\end{array} \right).
\end{eqnarray*}
Let
$y^{s,C}=[\mathbb{F}^{s,C}_{11},\mathbb{F}^{s,C}_{12};\mathbb{F}^{s,C}_{21},\mathbb{F}^{s,C}_{22}]$
be the predict output of CNN corresponding the input
$x^{s}=\mathbb{F}^{s}_{12}$, we can obtain the multi-static response
matrix as
\begin{eqnarray*} \mathbb{\widetilde{F}}^{s}_{full}=\left(
\begin{array}{cccc}\vspace{2mm}
\mathbb{F}^{s,C}_{11} & \mathbb{F}^{s}_{12}  \\
\mathbb{F}^{s,C}_{21} & \mathbb{F}^{s,C}_{22}
\end{array} \right).
\end{eqnarray*}
As a benchmark, we should compare the result with the nearest sample in the training set. Let
\begin{eqnarray*}
{s^*}=\mbox{arg}\min\limits_{1\leq n\leq n_1}
\left\{\|\mathbb{F}^{s}_{12}-\mathbb{F}_{12}^{n}\|^2_{F} \right\},
\end{eqnarray*}
for phased data, or
\begin{eqnarray*}
{s^*}=\mbox{arg}\min\limits_{1\leq n\leq n_1}
\left\{\| \left|\mathbb{F}^{s}_{12}\right| - \left|\mathbb{F}_{12}^{n} \right| \|^2_{F} \right\}
\end{eqnarray*}
for phaseless data. Then the nearest sample is given by
\begin{eqnarray*} \mathbb{\widetilde{F}}^{s^*}_{full}=\left(
\begin{array}{cccc}\vspace{2mm}
\mathbb{F}^{s^*}_{11} & \mathbb{F}^{s^*}_{12}  \\
\mathbb{F}^{s^*}_{21} & \mathbb{F}^{s^*}_{22}
\end{array} \right).
\end{eqnarray*}

\begin{figure}[!htbp]
\centering
\includegraphics[width=\textwidth]{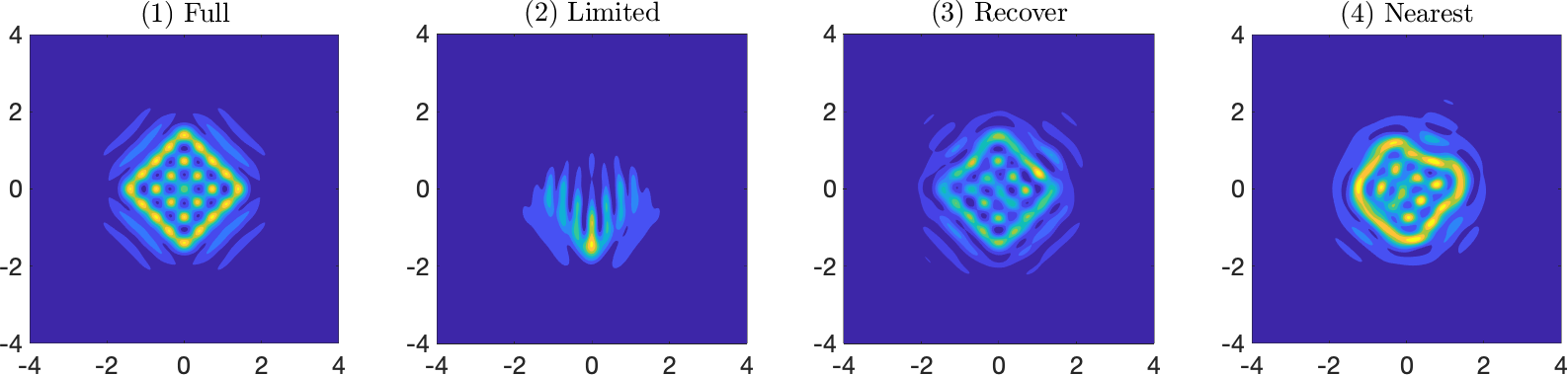}
\caption{Numerical constructions for the square-shaped obstacle by (1) the full
MSR matrix $\mathbb{F}^{s}_{full}$; (2) the limited MSR matrix
$\mathbb{F}^{s}_{12}$; (3) the recovery MSR matrix
$\mathbb{\widetilde{F}}^{s}_{full}$ via CNN; (4) the nearest sample
MSR matrix
$\mathbb{\widetilde{F}}^{s^*}_{full}$.}\label{fig:shape_square}
\end{figure}

Fig~\ref{fig:shape_square} shows the numerical constructions for
square-shaped obstacle by the full MSR matrix $\mathbb{F}^{s}_{full}$,
a limited MSR matrix $\mathbb{F}^{s}_{12}$, the recovery MSR matrix
$\mathbb{\widetilde{F}}^{s}_{full}$ via CNN, and the nearest sample
MSR matrix $\mathbb{\widetilde{F}}^{s^*}_{full}$. The corresponding
relative errors of recovery MSR matrix
$\mathbb{\widetilde{F}}^{s}_{full}$ and the nearest sample
MSR matrix $\mathbb{\widetilde{F}}^{s^*}_{full}$ are listed in Table~\ref{table:recover_near_ex2}.
Although the relative error of the test set $\bar{e}$ is around $13\%$, the numerical
construction by the recovery MSR matrix via CNN is very close to the numerical construction by the full MSR matrix, and much better than the reconstructions from the limited aperture data or the nearest sample.
\begin{table}[!htbp]
  \caption{The relative errors of recovery MSR matrix $\mathbb{\widetilde{F}}^{s}_{full}$ and the nearest sample MSR matrix $\mathbb{\widetilde{F}}^{s^*}_{full}$ in Example 2 with regularization term.} \label{table:recover_near_ex2}
\begin{tabular}{|c||c|c|c|}\hline
          & \multicolumn{3}{|c|}{relative error (\%) } \\ \hline
     matrix type      &    real   &    image       &  norm    \\ \hline \hline
   $\mathbb{\widetilde{F}}^{s}_{full}$    &   32.00   &    30.42       &  31.25   \\ \hline
   $\mathbb{\widetilde{F}}^{s^*}_{full}$    &   66.75   &    69.21       &  67.94  \\ \hline
\end{tabular}
\end{table}

\vspace{3mm}

%%%%%%%%%%%%%%%%%%%%%%%%%%%%%%%%%%%%%%%%%%

\noindent {\bf Example 3.} In this example, we consider the
phaseless obstacle reconstruction problem \eqref{eq:phaseless1}. This problem has less
information, and thus, is more difficult to reconstruct.

\begin{figure}[!htbp]
\centering
\includegraphics[width=\textwidth]{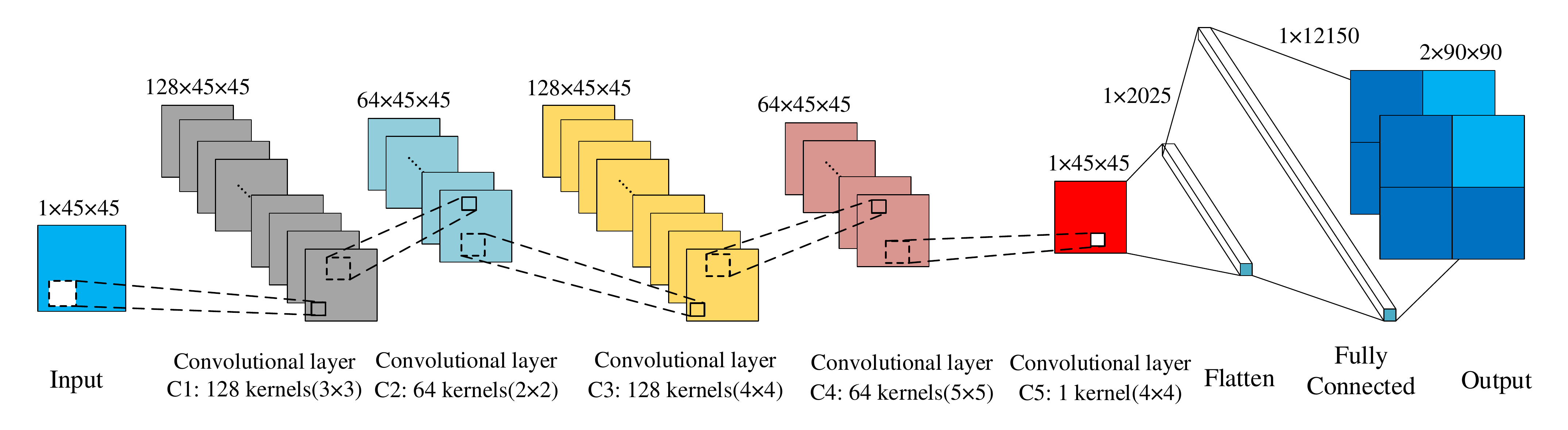}
\caption{The architecture of our convolution neural network in
Example 3.} \label{fig:EXAM3CNN}
\end{figure}

Here, the input data has just module information, that is, $x=\mathbb{F}^{n}_{12}\in \mathbb{R}^{m_1\times m_1}$, and the
target outputs remain unchanged
$y=[\mathbb{F}^{n}_{11},\mathbb{F}^{n}_{12};\mathbb{F}^{n}_{21},\mathbb{F}^{n}_{22}]\in
\mathbb{R}^{2m\times 2m \times 2}$. The architecture of convolution
neural network is similar to previous examples, except that the
input layer has only one matrix input (cf. Fig~\ref{fig:EXAM3CNN}).
The number of samples is $n^{*}=50000$, while the size of the train set
is $n_1=40000$ and the size of the test set is $n_2=10000$.

We adopt the same loss function $\mathcal{L}_2$ as in Example 2.
Table~\ref{table:phaseless_ex3} shows that the relative errors and $\overline{\mbox{MSE}}$s
of the test set are decreasing as $n_1$ and $n_2$ increase, while the
$\overline{\mbox{PSNR}}$s of the test set are increasing as $n_1$ and
$n_2$ increase. It is easy to see from Table~\ref{table:50T_ex2} and Table~\ref{table:BV_ex2} that,
the results of phaseless case are worse than that of phase case
because we have less information.
~\\
\begin{table}[!htbp]
  \caption{The relative errors, $\overline{\mbox{MSE}}$s, and $\overline{\mbox{PSNR}}$s of the test set in Example 3 with regularization term.} \label{table:phaseless_ex3}
\begin{tabular}{|r||c|c|c|c|c|c|c|c|c|}\hline
             & \multicolumn{3}{|c|}{relative error $\overline{e}$(\%) } &   \multicolumn{3}{|c|}{MSE($10^{-2}$)} &
    \multicolumn{3}{|c|}{PSNR}\\ \hline
   $(n_1,n_2)$    &    real   &    image       &  norm    &  real   & image  &  norm   & real    &image    & norm   \\ \hline \hline
   (5000,1250)    &   57.49   &    64.95       &  60.91   & 20.535   & 20.615  & 11.332   & 15.87   & 12.84   & 20.24  \\ \hline
   (10000,2500)   &   46.92   &    53.41       &  49.90   & 12.562   & 12.719  & 6.648   & 17.33   & 14.23   & 21.65  \\ \hline
   (20000,5000)   &   41.77   &    47.73       &  44.51   & 9.860   & 10.008  & 5.266   & 18.30   & 15.18   & 22.55  \\ \hline
   (40000,10000)  &   37.97   &    43.28       &  40.40   & 7.510   & 7.573  & 4.408   & 18.91   & 15.81   & 22.91  \\ \hline
\end{tabular}
\end{table}

As a test, we consider the unit circle centered at the origin, whose Fourier coefficients lies within the possible random numbers but not one of the \( 50000 \) samples in this example. Let $\mathbb{F}^{c}_{full}$, $\mathbb{F}^{c}_{12}$,
$\mathbb{\widetilde{F}}^{c}_{full}$,
$\mathbb{\widetilde{F}}^{c^*}_{full}$ be the corresponding full MSR
matrix, limited MSR matrix, recovery MSR matrix via CNN, and the
nearest sample MSR matrix, respectively, in the same manner as
in Example 2.
\begin{figure}[!htbp]
\centering
\includegraphics[width=\textwidth]{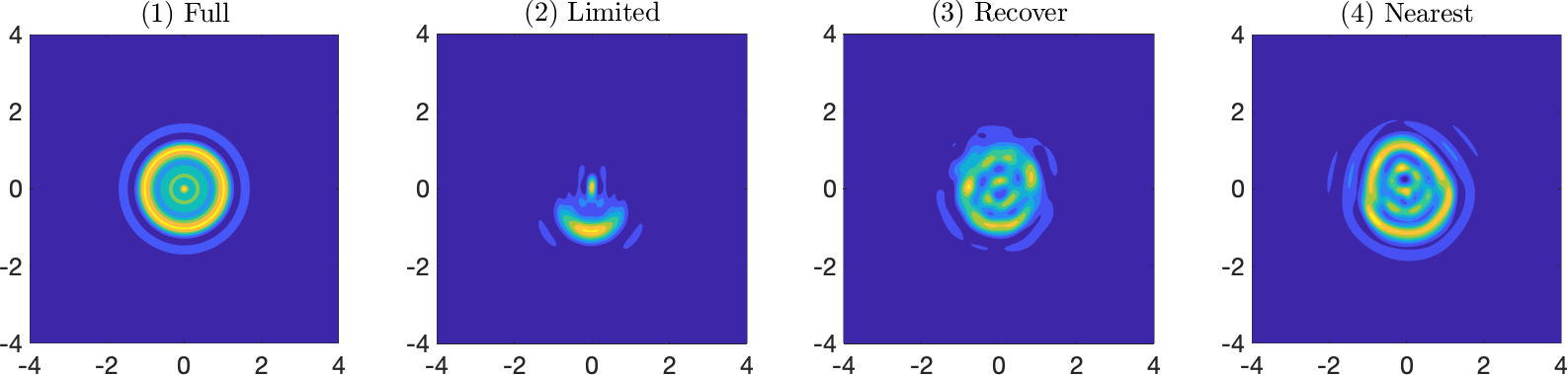}
\caption{Numerical constructions for circle-shape obstacle by (1) the full
MSR matrix $\mathbb{F}^{c}_{full}$; (2) the limited MSR matrix
$\mathbb{F}^{c}_{12}$; (3) the recovery MSR matrix
$\mathbb{\widetilde{F}}^{c}_{full}$ via CNN; (4) nearest sample
MSR matrix
$\mathbb{\widetilde{F}}^{c^*}_{full}$.}\label{fig:shape_circle}
\end{figure}
Fig~\ref{fig:shape_circle} shows the numerical
constructions of the circle by the full MSR matrix
$\mathbb{F}^{c}_{full}$, a limited MSR matrix $\mathbb{F}^{c}_{12}$,
the recovery MSR matrix $\mathbb{\widetilde{F}}^{c}_{full}$ via CNN, and
the nearest sample MSR matrix $\mathbb{\widetilde{F}}^{c^*}_{full}$.
Even for the phaseless case, the numerical construction by the recovery
MSR matrix via CNN is much better than the reconstruction from the original limited aperture data. The corresponding
relative errors of recovery MSR matrix
$\mathbb{\widetilde{F}}^{c}_{full}$ and the nearest sample
MSR matrix $\mathbb{\widetilde{F}}^{c^*}_{full}$ are listed in Table~\ref{table:recover_near_ex3}.
In fact, it is easy to see from Table~\ref{table:recover_near_ex3} that the recovered MSR matrix via CNN is slightly better than the nearest sample
recovery multi-static response matrix under the $\|\cdot\|_{F}$-norm.
\begin{table}[!htbp]
  \caption{The relative errors of recovery MSR matrix $\mathbb{\widetilde{F}}^{c}_{full}$ and the nearest sample MSR matrix $\mathbb{\widetilde{F}}^{c^*}_{full}$ in Example 3 with regularization term.} \label{table:recover_near_ex3}
\begin{tabular}{|c||c|c|c|}\hline
          & \multicolumn{3}{|c|}{relative error (\%) } \\ \hline
     matrix type      &    real   &    image       &  norm    \\ \hline \hline
   $\mathbb{\widetilde{F}}^{s}_{full}$    &   56.16   &    57.42       &  56.66   \\ \hline
   $\mathbb{\widetilde{F}}^{s^*}_{full}$    &   69.05   &    92.99       &  79.50  \\ \hline
\end{tabular}
\end{table}

\vspace{3mm}

%%%%%%%%%%%%%%%%%%%%%%%%%%%%%%%%%%%%%%%%%%

\noindent {\bf Example 4.} In this example, we consider the less
information of phase case. The architecture of convolution neural
network is similar to the one used in Example 1, except the input data
$x=\mathbb{F}^{n}_{12}\in \mathbb{R}^{m_1\times m_1 \times 2}$ with
$m_1=30$. The total number of samples is $n^{*}=50000$, while the
size of the train set is $n_1=40000$ and the size of the test set is
$n_2=10000$.

\begin{figure}[!htbp]
\centering
\includegraphics[width=\textwidth]{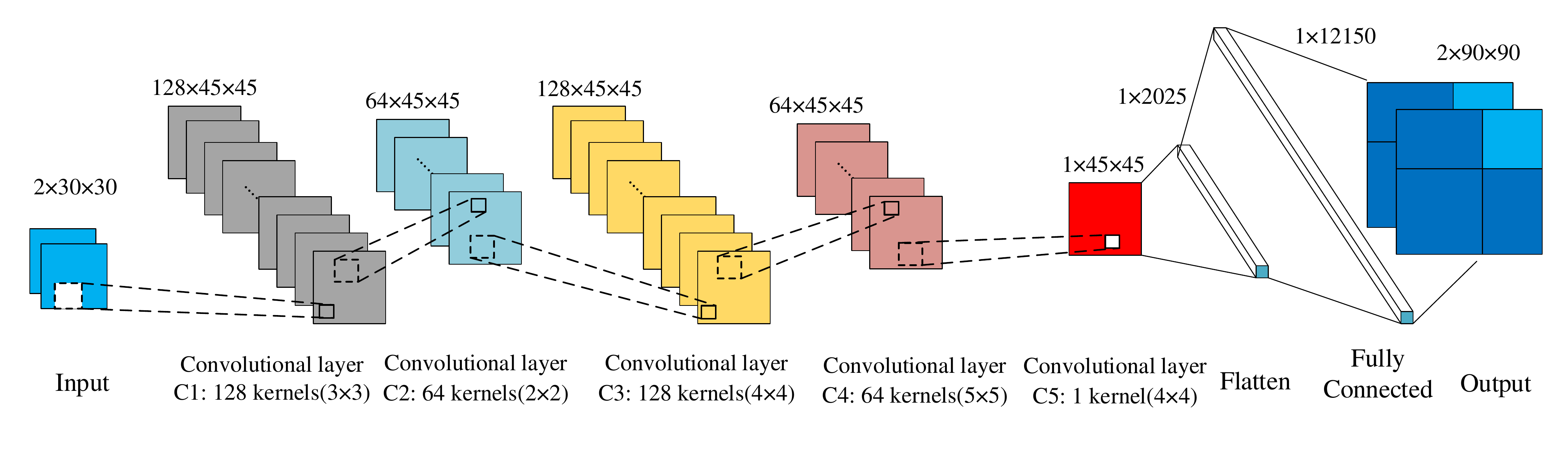}
\caption{The architecture of our convolution neural network in
Example 4.} \label{fig:EXAM4CNN}
\end{figure}

It is easy to see from Table~\ref{table:BV_ex2} and Table~\ref{table:phaseless_ex3} that, the results of
$m_1=30$ are worse than that of $m_1=45$ because that we have less
information.

\begin{table}[!htbp]
  \caption{The relative errors, $\overline{\mbox{MSE}}$s, and $\overline{\mbox{PSNR}}$s of test set in Example 4 with regularization term and input data $m_1=30$.} \label{table:input30_ex4}
\begin{tabular}{|r||c|c|c|c|c|c|c|c|c|}\hline
    error         & \multicolumn{3}{|c|}{relative error $\overline{e}$(\%) } &   \multicolumn{3}{|c|}{$\overline{\mbox{MSE}}$($10^{-2}$)} &
    \multicolumn{3}{|c|}{$\overline{\mbox{PSNR}}$}\\ \hline
   $(n_1,n_2)$    &    real   &    image       &  norm    &  real   & image  &  norm   & real    &image    & norm   \\ \hline \hline
   (5000,1250)    &   30.72   &    34.75       &  32.58   & 5.572   & 5.657  & 3.821   & 21.17   & 18.18   & 23.99  \\ \hline
   (10000,2500)   &   24.86   &    28.13       &  26.37   & 3.505   & 3.525  & 2.542   & 22.81   & 19.80   & 25.49  \\ \hline
   (20000,5000)   &   18.93   &    21.37       &  20.06   & 2.116   & 2.128  & 1.652   & 25.25   & 22.27   & 27.54  \\ \hline
   (40000,10000)  &   15.14   &    17.11       &  16.05   & 1.295   & 1.290  & 1.029   & 27.05   & 24.04   & 29.41  \\ \hline
\end{tabular}
\end{table}

\begin{figure}[!htbp]
\centering
\includegraphics[width=.6\textwidth]{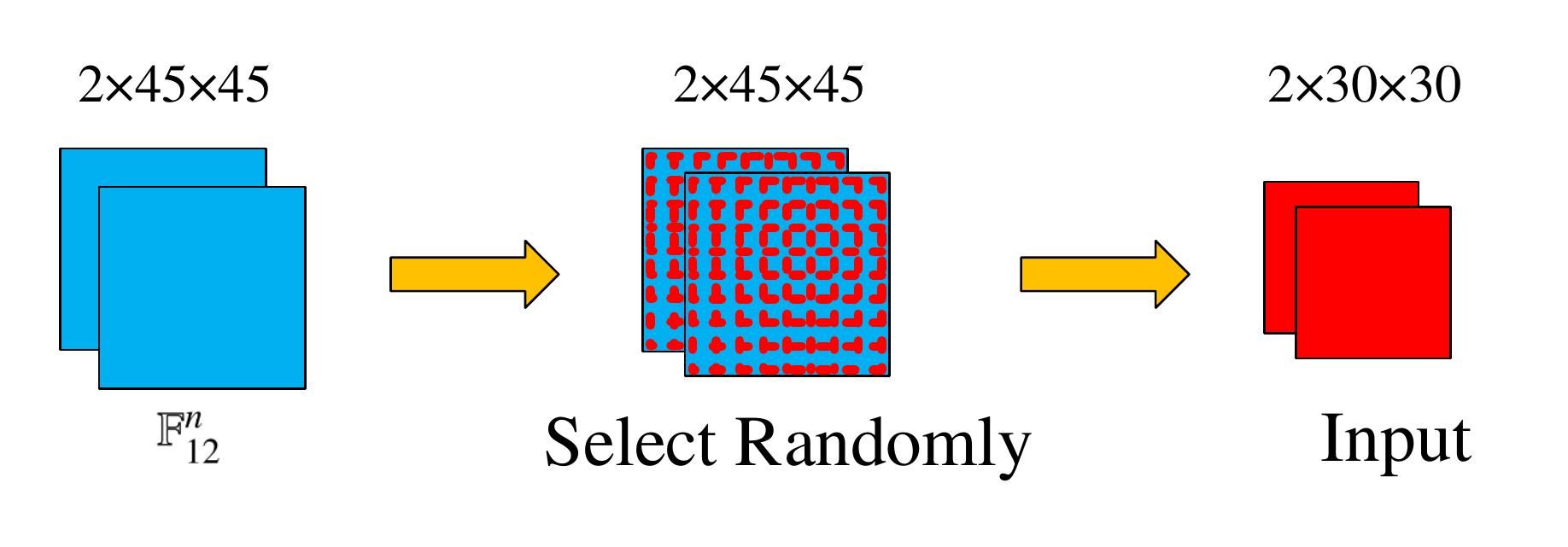}
\caption{The architecture of our convolution neural network in
Example 4 with input being chosen randomly.} \label{fig:EXAM5CNN}
\end{figure}

In order to show the robustness of our algorithm, we choose $30$
rows and $30$ columns randomly in the sub-block $\mathbb{F}_{12}\in
C^{m\times m}$ with $m=45$. The architecture of convolution neural
network is similar to Fig~\ref{fig:EXAM4CNN}, except the input data
$x=\mathbb{F}^{n}_{12}\in \mathbb{R}^{m_1\times m_1 \times 2}$ with
$m_1=30$ chosen randomly (see Fig~\ref{fig:EXAM5CNN}).

It is easy to see from Table~\ref{table:BV_ex2} and Table~\ref{table:input30_ex4} that, the results of
using randomly chosen inputs with $m_1=30$ are worse than that of $m_1=45$ because
we have less information. But the results of $m_1=30$ chosen
randomly are only slightly different from that of using static inputs as $n_1$
and $n_2$ increase (Table~\ref{table:input30_ex4} and Table~\ref{table:random30_ex4}).
\begin{table}[!htbp]
  \caption{The relative errors, $\overline{\mbox{MSE}}$s, and $\overline{\mbox{PSNR}}$s of test set in Example 4 with regularization term and input data $m_1=30$ chosen randomly.}\label{table:random30_ex4}
\begin{tabular}{|r||c|c|c|c|c|c|c|c|c|}\hline
    error         & \multicolumn{3}{|c|}{relative error $\overline{e}$(\%) } &   \multicolumn{3}{|c|}{$\overline{\mbox{MSE}}$($10^{-2}$)} &
    \multicolumn{3}{|c|}{$\overline{\mbox{PSNR}}$}\\ \hline
   $(n_1,n_2)$    &    real   &    image       &  norm    &  real   & image  &  norm   & real    &image    & norm   \\ \hline \hline
   (5000,1250)    &   33.79   &    37.90       &  35.70   & 6.580   & 6.541  & 4.152   & 20.31   & 17.39   & 23.32  \\ \hline
   (10000,2500)   &   27.44   &    31.22       &  29.20   & 4.295   & 4.371  & 3.020   & 22.01   & 18.95   & 24.67  \\ \hline
   (20000,5000)   &   20.36   &    23.06       &  21.62   & 2.495   & 2.509  & 1.936   & 24.69   & 21.67   & 27.01  \\ \hline
   (40000,10000)  &   15.24   &    17.24       &  16.16   & 1.340   & 1.334  & 1.115   & 27.07   & 24.05   & 29.27  \\ \hline
\end{tabular}
\end{table}

\section*{Acknowledgment}
The work of K. Zhang was supported by
the NSF of China under the grant No. 11871245, 11771179, and by the
Key Laboratory of Symbolic Computation and Knowledge Engineering of
Ministry of Education, Jilin University (93K172018Z01). The authors would like to acknowledge the helpful suggestion and discussion with Prof.
Hongyu Liu of CityU, Prof. Yuliang Wang of HKBU, and Prof. Xiaodong Liu of AMT.

\end{document}